\documentclass[preprint,amsmath,amssymb,aps,showpacs]{revtex4-1} 
\usepackage{graphicx}
\usepackage{amssymb}
\usepackage{amsmath}
\usepackage{amsbsy}
\usepackage{bm}

\begin{document}
\setcounter{page}{1}
\newtheorem{t1}{Theorem}[section]
\newtheorem{d1}{Definition}[section]
\newtheorem{c1}{Corollary}[section]
\newtheorem{l1}{Lemma}[section]
\newtheorem{r1}{Remark}[section]

\newcommand{\cA}{{\cal A}}
\newcommand{\cB}{{\cal B}}
\newcommand{\cC}{{\cal C}}
\newcommand{\cD}{{\cal D}}
\newcommand{\cE}{{\cal E}}
\newcommand{\cF}{{\cal F}}
\newcommand{\cG}{{\cal G}}
\newcommand{\cH}{{\cal H}}
\newcommand{\cI}{{\cal I}}
\newcommand{\cJ}{{\cal J}}
\newcommand{\cK}{{\cal K}}
\newcommand{\cL}{{\cal L}}
\newcommand{\cM}{{\cal M}}
\newcommand{\cN}{{\cal N}}
\newcommand{\cO}{{\cal O}}
\newcommand{\cP}{{\cal P}}
\newcommand{\cQ}{{\cal Q}}
\newcommand{\cR}{{\cal R}}
\newcommand{\cS}{{\cal S}}
\newcommand{\cT}{{\cal T}}
\newcommand{\cU}{{\cal U}}
\newcommand{\cV}{{\cal V}}
\newcommand{\cX}{{\cal X}}
\newcommand{\cW}{{\cal W}}
\newcommand{\cY}{{\cal Y}}
\newcommand{\cZ}{{\cal Z}}

\def\cl{\centerline}
\def\bd{\begin{description}}
\def\be{\begin{enumerate}}
\def\ben{\begin{equation}}

\def\een{\end{equation}}

\def\benr{\begin{eqnarray}}
\def\eenr{\end{eqnarray}}
\def\benrr{\begin{eqnarray*}}
\def\eenrr{\end{eqnarray*}}
\def\bes {\begin{equation*}}
\def\ees {\end{equation*}}
\def\bbm {\begin{bmatrix}}
\def\ebm {\end{bmatrix}}
\def\ed{\end{description}}
\def\ee{\end{enumerate}}
\def\al{\alpha}
\def\b{\beta}
\def\bR{\bar\R}
\def\bc{\begin{center}}
\def\ec{\end{center}}
\def\d{\dot}
\def\D{\Delta}
\def\del{\delta}
\def\ep{\epsilon}
\def\g{\gamma}
\def\G{\Gamma}
\def\h{\hat}
\def\iny{\infty}
\def\La{\Longrightarrow}
\def\la{\lambda}
\def\m{\mu}
\def\n{\nu}
\def\noi{\noindent}
\def\Om{\Omega}
\def\om{\omega}
\def\p{\psi}
\def\pr{\prime}
\def\r{\ref}
\def\R{{\bf R}}
\def\ra{\rightarrow}
\def\s{\sum_{i=1}^n}
\def\si{\sigma}
\def\Si{\Sigma}
\def\t{\tau}
\def\th{\theta}
\def\Th{\Theta}

\def\vep{\varepsilon}
\def\vp{\varphi}
\def\pa{\partial}
\def\un{\underline}
\def\ov{\overline}
\def\fr{\frac}
\def\sq{\sqrt}

\def\WW{\begin{stack}{\circle \\ W}\end{stack}}
\def\ww{\begin{stack}{\circle \\ w}\end{stack}}
\def\st{\stackrel}
\def\Ra{\Rightarrow}
\def\R{{\mathbb R}}
\def\bi{\begin{itemize}}
\def\ei{\end{itemize}}
\def\i{\item}
\def\bt{\begin{tabular}}
\def\et{\end{tabular}}
\def\lf{\leftarrow}
\def\nn{\nonumber}
\def\va{\vartheta}
\def\wh{\widehat}
\def\vs{\vspace}
\def\Lam{\Lambda}
\def\sm{\setminus}
\def\ba{\begin{array}}
\def\ea{\end{array}}
\def\bd{\begin{description}}
\def\ed{\end{description}}
\def\lan{\langle}
\def\ran{\rangle}
\def\mf{\mathbf}
\def\ts{\times}
\def\ots{\otimes}
\def\bs{\boldsymbol}
\def\l{\label}
\def\r{\ref}
\numberwithin{equation}{section}

\title{On Entanglement and Separability}

\author{Dhananjay P. Mehendale }
\email{dhananjay.p.mehendale@gmail.com}
\affiliation{Department of Physics, Savitribai Phule University
of Pune,Pune,India-411007.}

\date{\today}

\begin{abstract}
  
We present a necessary and sufficient condition to determine the entanglement status of an arbitrary $N$-qubit 
quantum state (may be pure or mixed) represented by the density matrix, $\rho_N.$ We develop a new approach and a new 
criterion for the problem of deciding entanglement status. 
Further, we develop as an important application of entanglement a new quantum protocol for superluminal communication 
of classical information in terms of a desired ordered sequence of classical bits. 
We then show that this new quantum protocol for superluminal communication of classical information can be used in the 
quantum teleportation protocol \cite{CB} for achieving superluminal quantum teleportation.

\end{abstract}
 

\maketitle

\section{Introduction}

One of the central issues in quantum information theory is whether a given 
$N$-qubit pure or mixed quantum state represented by the density matrix, $\rho_N,$ is separable or 
entangled \cite{Peres1,Horo,PHoro,Horos,wer,alb,bou,nie,per}. 
This important question of deciding whether a given 
$N$-qubit (pure or mixed) quantum state represented by the density matrix, $\rho_N,$ is separable 
or entangled is completely solved in this paper. In this paper, we present a criterion which leads 
to an algorithm for deciding the entanglement status of an arbitrary $N$-qubit pure or mixed quantum state, 
represented in terms of a density matrix, $\rho_N,$ of size $2^N\ts 2^N.$ 

Our result for $N$-qubit ``pure'' quantum state is as follows: A quantum state represented by a density 
matrix, $\rho_N$ of size $2^N\ts 2^N$ is a ``pure separable'' state if and only if there is a matrix called 
{\bf basic matrix} of size $4\ts 4$ whose rank is equal to one and other matrices of size $4\ts 4$ 
(in all $4^{(N-2)}$ in number including the basic matrix) made up of certain coefficients that arise in the 
representation for this density matrix, $\rho_N,$ in terms of the Generalized Pauli Basis and all these  
other matrices are certain constant multiples of this {\bf basic matrix}.  

Our result for $N$-qubit ``mixed'' quantum state is as follows: A quantum state 
represented by a density matrix, $\rho_N$ of size $2^N\ts 2^N$ is a ``mixed separable'' state 
if and only if there is a matrix called {\bf basic matrix} of size $4\ts 4$ and other matrices of size $4\ts 4$ 
(in all $4^{(N-2)}$ in number including the basic matrix) made up of certain coefficients that arise in the 
representation for this density matrix, $\rho_N,$ in terms of the Generalized Pauli Basis such that all the other
matrices are certain constant multiples of this {\bf basic matrix} and this {\bf basic matrix} can be split into 
sum of certain other matrices of rank one and of size $4\ts 4.$  

Entanglement describes a correlation between different parts of a quantum
system that exceeds anything that is possible classically. When a quantum system is
in such an entangled state, actions performed on one sub-system will have 
effects on another sub-system even though that sub-system is not acted upon directly and could be far away. 
This leads to many highly counterintuitive phenomena. All the known quantum 
algorithms that display an exponential speedup over their
classical counterparts exploit such entanglement-induced side effects in one way or
another \cite{Colin}. Therefore, to study, understand, and characterize entanglement is one of the very important problems in 
quantum information theory. 
We also will give in this paper a new quantum protocol for instantaneous and simultaneous transmission of classical information to faraway 
locations. This new quantum protocol achieves instantaneous and simultaneous transfer of the ``same''  
classical information in terms of certain desired sequence of classical bits from the data transmission center 
to different faraway locations from the data transmission center. Further, these locations are not only faraway 
from data transmission center but also could be located faraway 
from one another too. This new quantum protocol thus serves as a superluminal communicator for sending any desired 
classical information. This classical information to be transmitted is not known a priory to the data transmission center and 
it arrives there from time to time for transmission purpose when Alice at data transmission center and other participants 
have already gone faraway from one another. 
We further show that the instantaneous part of the ability to transmit the classical information in terms of a sequence 
of classical bits also enables one to achieve instantaneous teleportation of a quantum state 
to some (single) chosen faraway location which was impossible for the teleportation protocol \cite{CB} because there was no 
quantum protocol available to achieve superluminal communication of classical information in terms of classical bits. The 
two classical bits that appear after Alice's Bell basis measurement of her qubits are needed by Bob, now 
stationed at some faraway location, to complete the teleportation protocol \cite{CB}. Our new quantum protocol for 
superluminal classical communication makes possible the needed instantaneous transmission of these two classical bits 
to Bob to successfully complete the teleportation protocol \cite{CB}. 

\section{Generalized Pauli Basis}

Let $\{I(= \sigma_0),\sigma_1,\sigma_2,\sigma_3\}$ denote the well-known Pauli matrices:   
\bes
\ I
\ =
\bbm
1  &  0 \\
0  &  1
\ebm
\ees
\bes
\ \sigma_1
\ =
\bbm
0  &  1 \\
1  &  0
\ebm
\ees
\bes
\ \sigma_2
\ =
\bbm
0  &  -i \\
i  &  0
\ebm
\ees
\bes
\ \sigma_3
\ =
\bbm
1  &  0 \\
0  &  -1
\ebm
\ees
Let $A=[a_{ij}]$ and $B=[b_{ij}]$ be some $2\ts 2$ hermitian matrices then we denote their scalar or dot product as $A.B$ and define it as  
$A.B = \sum_{i=1}^2\sum_{j=1}^2 a_{ij}^*b_{ij},$ where $a_{ij}^*$ denotes the complex conjugate of $a_{ij}.$ 
Pauli matrices $\{I (= \sigma_0),\sigma_1,\sigma_2,\sigma_3\}$ form the orthogonal vector 
space basis for $2\ts 2$ hermitian matrices i.e. every  $2\ts 2$ hermitian matrix can be uniquely expressed 
as linear combination of Pauli matrices over reals, $\mathbb R.$ Thus, if $A = \sum_{i=1}^4\alpha_i\sigma_i$ then 
$\alpha_i = A.\sigma_i/\sigma_i.\sigma_i$ It is easy to check that $\sigma_i.\sigma_i = 2$ for all $i,$ 
and $\sigma_i.\sigma_j = 0$ if $i\neq j$ $i,j\in \{0,1,2,3\}.$ 
The density operator is a positive semidefinite operator with nonnegative eigenvalues.
The density matrix, $\rho_N,$ corresponding to an $N$-qubit state is a $2^N\ts 2^N$ matrix with certain well defined 
properties: For example, (i) the sum of the diagonal elements of $\rho,$ tr$(\rho),$ is equal to one. (ii) $\rho$ is hermitian, (iii) If $\rho$ corresponds 
to a pure state then tr$(\rho^2) = 1,$ (iv) If $\rho$ corresponds 
to a mixed state then tr$(\rho^2) < 1,$ etc. 

We now define the Generalized Pauli matrices.

{\bf Definition 2.1:} Generalized Pauli matrices of size $2^N\ts 2^N$ are those obtained from usual 
Pauli matrices, $\{I (= \sigma_0),\sigma_1,\sigma_2,\sigma_3\},$ by taking their all possible tensor products: 
$G_{i_1i_2\cdots i_s\cdots i_N} = \sigma_{i_1}\ots \sigma_{i_2}\ots \cdots\ots \sigma_{i_s}\ots \cdots\ots \sigma_{i_N},$ 
where each $\sigma_{i_s}$ is a usual Pauli 
matrix i.e. $\sigma_{i_s} \in \{I(= \sigma_0),\sigma_1,\sigma_2,\sigma_3\}.$ 

Now, let $C=[c_{ij}]$ and $D=[d_{ij}]$ be some $2^N\ts 2^N$ hermitian matrices then we denote their scalar or dot product as $C.D$ and define it as  
$C.D = \sum_{i=1}^{2^N}\sum_{j=1}^{2^N} c_{ij}^*d_{ij},$ where $c_{ij}^*$ denotes the complex conjugate of $c_{ij}.$ 
Similar to the well-known simple result stated 
above, namely, Pauli matrices $\{I (= \sigma_0),\sigma_1,\sigma_2,\sigma_3\}$ form the orthogonal vector 
space basis for $2\ts 2$ hermitian matrices i.e. every  $2\ts 2$ hermitian matrix can be uniquely expressed 
as linear combination of Pauli matrices over reals, $\mathbb R,$ one can easily obtain the following generalization 
by proceeding on similar lines: Generalized Pauli 
matrices $G_{i_1i_2\cdots i_r\cdots i_N} = \sigma_{i_1}\ots \sigma_{i_2}\ots \cdots\ots \sigma_{i_r}\ots \cdots\ots \sigma_{i_N},$ 
where each $\sigma_{i_r}$ is a usual Pauli 
matrix i.e. $\sigma_{i_r} \in \{I(= \sigma_0),\sigma_1,\sigma_2,\sigma_3\},$ form the orthogonal vector 
space basis for $2^N\ts 2^N$ hermitian matrices i.e. every  $2^N\ts 2^N$ hermitian matrix can be uniquely expressed 
as linear combination of Generalized Pauli matrices over reals, $\mathbb R.$ Let $A$ be $2^N\ts 2^N$ hermitian matrix 
then $A$ can be uniquely expressed as 
$$A = \sum_{i_1,i_2,\cdots,i_r,\cdots,i_N} \alpha_{i_1i_2\cdots i_r\cdots i_N}G_{i_1i_2\cdots i_r\cdots i_N}$$
where 
$$\alpha_{i_1i_2\cdots i_r\cdots i_N} = (A).(G_{i_1i_2\cdots i_r\cdots i_N})/(G_{i_1i_2\cdots i_r\cdots i_N}).(G_{i_1i_2\cdots i_r\cdots i_N})$$
It is easy to check that the scalar or dot product  
$(G_{i_1i_2\cdots i_r\cdots i_N}).(G_{i_1i_2\cdots i_r\cdots i_N}) = 2^N$ and also the scalar or dot product 
$(G_{i_1i_2\cdots i_r\cdots i_N}).(G_{j_1j_2\cdots j_r\cdots j_N}) = 0$ if some $i_r \neq j_r,$ $r\in \{1,2,\cdots, ,N\}.$ 

The vectors in the Generalized Pauli basis for $N$-dimensional case are actually the matrices of size $2^N\ts 2^N$ and 
they have a special form. 

(i) For the case of $N=2$ these ``sixteen'' matrices (each of these matrices are of size $4\ts 4$ and have 
exactly ``four'' nonzero elements $\in \mathbb{C},$ the field of complex numbers) 
representing basis vectors have ``four'' types of forms and so 
they can be divided into ``four'' (independent and disjoint) groups as given below: 

{\bf Group 1:} The $4\ts 4$ matrices representing basis vectors in

$(I,\sigma_3)\ots (I,\sigma_3) = I\ots I, I\ots \sigma_3, \sigma_3\ots I, \sigma_3\ots \sigma_3$ 

and they have the form of a ``Diagonal'' matrix:
\bes
\bbm
a_{11} &  0 & 0 & 0  \\
0  &  a_{22} & 0 & 0  \\
0  &  0 & a_{33} & 0  \\
0  &  0 & 0 & a_{44} 
\ebm
\ees 
{\bf Group 2:} The $4\ts 4$ matrices representing basis vectors in 

$(I,\sigma_3)\ots (\sigma_1,\sigma_2) = I\ots \sigma_1, I\ots \sigma_2, \sigma_3\ots \sigma_1, \sigma_3\ots \sigma_2$

and they have the form as given below: 
\bes
\bbm
0  &  a_{12} & 0 & 0  \\
a_{21}  & 0 & 0 & 0  \\
0  & 0 & 0 & a_{34}  \\
0  &  0 & a_{43} & 0 
\ebm
\ees
{\bf Group 3:} The $4\ts 4$ matrices representing basis vectors in

$(\sigma_1,\sigma_2)\ots (I,\sigma_3) = \sigma_1\ots I, \sigma_1\ots \sigma_3, \sigma_2\ots I, \sigma_2\ots \sigma_3$

and they have the form as given below: 
\bes
\bbm
0  &  0 & a_{13} & 0  \\
0  &  0 & 0 & a_{24}  \\
a_{31}  &  0 & 0 & 0  \\
0  &  a_{42} & 0 & 0 
\ebm
\ees
{\bf Group 4:} The $4\ts 4$ matrices representing basis vectors in

$(\sigma_1,\sigma_2)\ots (\sigma_1,\sigma_2) = \sigma_1\ots \sigma_1\, \sigma_1\ots \sigma_2, \sigma_2\ots \sigma_1\, \sigma_2\ots \sigma_2$

and they have the form of an ``Off-Diagonal'' matrix:
\bes
\bbm
0  &  0 & 0 & a_{14}  \\
0  &  0 & a_{23} & 0  \\
0  &  a_{32} & 0 & 0  \\
a_{41}  &  0 & 0 & 0 
\ebm
\ees
(ii) For the case of $N=3$ these ``sixtyfour'' matrices (each of these matrices are of size $8\ts 8$ and have 
exactly ``eight'' nonzero elements $\in \mathbb{C},$ the field of complex numbers) representing basis vectors have ``eight'' types of forms and so 
they can be divided into ``eight'' (independent and disjoint) groups as given below: 

{\bf Group 1:} The $8\ts 8$ matrices representing basis vectors in

$(I,\sigma_3)\ots (I,\sigma_3)\ots (I,\sigma_3),$ are ``eight'' in number 

and they have the form of a ``Diagonal'' matrix.

{\bf Group 2:} The $8\ts 8$ matrices representing basis vectors in  

$(I,\sigma_3)\ots (I,\sigma_3)\ots (\sigma_1,\sigma_2),$ are ``eight'' in number 

and they have the form as given below: 
\bes
\bbm
0  &  a_{12} & 0 & 0 & 0 & 0 & 0 & 0 \\
a_{21}  &  0 & 0 & 0 & 0 & 0 & 0 & 0 \\
0  &  0 & 0 & a_{34} & 0 & 0 & 0 & 0 \\
0  &  0 & a_{43} & 0 & 0 & 0 & 0 & 0 \\
0  &  0 & 0 & 0 & 0 & a_{56} & 0 & 0 \\
0  &  0 & 0 & 0 & a_{65} & 0 & 0 & 0 \\
0  &  0 & 0 & 0 & 0 & 0 & 0 & a_{78} \\
0  &  0 & 0 & 0 & 0 & 0 & a_{87} & 0 
\ebm
\ees
{\bf Group 3:} The $8\ts 8$ matrices representing basis vectors in 

$(I,\sigma_3)\ots (\sigma_1,\sigma_2)\ots (I,\sigma_3),$ are ``eight'' in number 

and they have the form as given below: 
\bes
\bbm
0  &  0 & a_{13} & 0 & 0 & 0 & 0 & 0 \\
0  &  0 & 0 & a_{24} & 0 & 0 & 0 & 0 \\
a_{31}  &  0 & 0 & 0 & 0 & 0 & 0 & 0 \\
0  &  a_{42} & 0 & 0 & 0 & 0 & 0 & 0 \\
0  &  0 & 0 & 0 & 0 & 0 & a_{57} & 0 \\
0  &  0 & 0 & 0 & 0 & 0 & 0 & a_{68} \\
0  &  0 & 0 & 0 & a_{75} & 0 & 0 & 0 \\
0  &  0 & 0 & 0 & 0 & a_{86} & 0 & 0 
\ebm
\ees
{\bf Group 4:} The $8\ts 8$ matrices representing basis vectors in 

$(I,\sigma_3)\ots (\sigma_1,\sigma_2)\ots (\sigma_1,\sigma_2),$ are ``eight'' in number 

and they have the form as given below: 
\bes
\bbm
0 & 0 & 0 & a_{14} & 0 & 0 & 0 & 0 \\
0 & 0 & a_{23} & 0 & 0 & 0 & 0 & 0 \\
0 & a_{32} & 0 & 0 & 0 & 0 & 0 & 0 \\
a_{41} & 0 & 0 & 0 & 0 & 0 & 0 & 0 \\
0 & 0 & 0 & 0 & 0 & 0 & 0 & a_{58} \\
0 & 0 & 0 & 0 & 0 & 0 & a_{67} & 0 \\
0 & 0 & 0 & 0 & 0 & a_{76} & 0 & 0 \\
0 & 0 & 0 & 0 & a_{85} & 0 & 0 & 0 
\ebm
\ees
{\bf Group 5:} The $8\ts 8$ matrices representing basis vectors in 

$(\sigma_1,\sigma_2)\ots (I,\sigma_3)\ots (I,\sigma_3),$ are ``eight'' in number 

and they have the form as given below:  
\bes
\bbm
0 & 0 & 0 & 0 & a_{15} & 0 & 0 & 0 \\
0 & 0 & 0 & 0 & 0 & a_{26} & 0 & 0 \\
0 & 0 & 0 & 0 & 0 & 0 & a_{37} & 0 \\
0 & 0 & 0 & 0 & 0 & 0 & 0 & a_{48} \\
a_{51} & 0 & 0 & 0 & 0 & 0 & 0 & 0 \\
0 & a_{62} & 0 & 0 & 0 & 0 & 0 & 0 \\
0 & 0 & a_{73} & 0 & 0 & 0 & 0 & 0 \\
0 & 0 & 0 & a_{84} & 0 & 0 & 0 & 0 
\ebm
\ees
{\bf Group 6:} The $8\ts 8$ matrices representing basis vectors in 

$(\sigma_1,\sigma_2)\ots (I,\sigma_3)\ots (\sigma_1,\sigma_2),$ are ``eight'' in number 

and they have the form as given below: 
\bes
\bbm
0 & 0 & 0 & 0 & 0 & a_{16} & 0 & 0 \\
0 & 0 & 0 & 0 & a_{25} & 0 & 0 & 0 \\
0 & 0 & 0 & 0 & 0 & 0 & 0 & a_{38} \\
0 & 0 & 0 & 0 & 0 & 0 & a_{47} & 0 \\
0 & a_{52} & 0 & 0 & 0 & 0 & 0 & 0 \\
a_{61} & 0 & 0 & 0 & 0 & 0 & 0 & 0 \\
0 & 0 & 0 & a_{74} & 0 & 0 & 0 & 0 \\
0 & 0 & a_{83} & 0 & 0 & 0 & 0 & 0 
\ebm
\ees
{\bf Group 7:} The $8\ts 8$ matrices representing basis vectors in 

$(\sigma_1,\sigma_2)\ots (\sigma_1,\sigma_2)\ots (I,\sigma_3),$ are ``eight'' in number 

and they have the form as given below: 
\bes
\bbm
0 & 0 & 0 & 0 & 0 & 0 & a_{17} & 0 \\
0 & 0 & 0 & 0 & 0 & 0 & 0 & a_{28} \\
0 & 0 & 0 & 0 & a_{35} & 0 & 0 & 0 \\
0 & 0 & 0 & 0 & 0 & a_{46} & 0 & 0 \\
0 & 0 & a_{53} & 0 & 0 & 0 & 0 & 0 \\
0 & 0 & 0 & a_{64} & 0 & 0 & 0 & 0 \\
a_{71} & 0 & 0 & 0 & 0 & 0 & 0 & 0 \\
0 & a_{82} & 0 & 0 & 0 & 0 & 0 & 0 
\ebm
\ees
{\bf Group 8:} The $8\ts 8$ matrices representing basis vectors in  

$(\sigma_1,\sigma_2)\ots (\sigma_1,\sigma_2)\ots (\sigma_1,\sigma_2),$ are ``eight'' in number 

and they have the form of an ``Off-Diagonal'' matrix.

\section{Preliminaries}

All 1-qubit states represented by their corresponding $2\ts 2$ density matrices being non-composite 
are separable. On the other hand $N$-qubit ($N \ge 2$) states represented by their 
corresponding $2^N\ts 2^N$ density matrices 
being composite can be separable or entangled. An $N$-qubit register is in a state of the form:
$$|\psi\ran = \sum_{i_1,i_2,\cdots, i_N \in \{0,1\}} c_{i_1i_2\cdots i_N} |i_1i_2\cdots i_N\ran$$ 
such that $\sum |c_{i_1i_2\cdots i_N}|^2 =1.$ This state expressible as a superposition of orthonormal basis states 
is said to be a ``pure'' state. We have complete knowledge of this $N$-qubit state. The states about which we do not 
have complete knowledge and which are statistical mixture of certain pure states and further these pure states are 
not necessarily mutually orthogonal are called ``mixed'' states. Thus, we have only partial (probabilistic) knowledge about 
these mixed states and we only know that a mixed state is in one of the (not necessarily orthogonal)
states among the states $|\psi_1\ran,|\psi_2\ran,\cdots, |\psi_N\ran$ with probabilities $p_1,p_2,\cdots, p_N$ respectively
such that $\sum p_i =1.$
The mixed states being a statistical mixture of pure states we are therefore a little uncertain of what a mixed  
state actually is. The best description of such a mixed state is given by the density operator: 
$$\rho_N = \sum_{i=1}^N p_i|\psi_i\ran\lan \psi_i|$$ where $\sum p_i =1$ and in general 
$\lan \psi_i|\psi_j\ran \neq 0$ when $i\neq j.$ Note that the decomposition of a given density operator into a weighted 
sum of pure states is non-unique. Any decomposition that synthesizes the density operator is as legitimate
as any other. This means that there is no unique mixed state to which each density operator corresponds. 

Entanglement describes a correlation between different parts of a quantum system which persists regardless of how far 
apart these parts scatter. Some tasks such as teleportation of a quantum state use entanglement in an essential way. 
It is therefore important to develop simple and useful criteria to determine which states are entangled and which are not. 
The states which are not entangled are called separable states. 

{\bf Separable States:} 

{\bf Definition 3.1:} A pure (mixed) bipartite state, represented by the wavefunction, $|\psi_{AB}\ran,$ (represented by the density operator, $\rho_{AB},$) of a composite quantum system defined on a 
Hilbert space $H_A\ots H_B$ is separable if and only if it can be written as $|\psi_{AB}\ran = |\psi_A\ran\ots |\psi_B\ran$ 
(it can be represented or approximated by $\rho_{AB} = \sum p_i \rho_i^A\ots\rho_i^B$ where $\sum p_i =1$). 

{\bf Definition 3.2:} A pure (mixed) multipartite state, represented by the wavefunction, $|\psi_{A_1A_2\cdots A_N}\ran,$ (represented by the density operator, $\rho_{A_1A_2\cdots A_N}$) of a composite 
quantum system defined on a Hilbert space $H_{A_1}\ots H_{A_2}\ots \cdots\ots H_{A_N}$ is separable if and only if it 
can be written as 
$|\psi_{A_1A_2\cdots A_N}\ran = |\psi_{A_1}\ran\ots |\psi_{A_2}\ran\ots \cdots\ots |\psi_{A_N}\ran$ 
(it can be represented or approximated by 
$\rho_{A_1A_2\cdots A_N} = \sum p_i \rho_i^{A_1}\ots\rho_i^{A_2}\ots \cdots\ots \rho_i^{A_N}$ where $\sum p_i =1$). 

{\bf Entangled States:} 

If a state of a composite quantum system is not a separable state then it is an entangled state.

States can be of four possible types: Separable Pure, Separable Mixed, Entangled Pure, Entangled Mixed. We make use of the 
above given definitions and give a constructive procedure to check the type of the given state. 

The density matrix corresponding to 1-qubit state is a $2^1\ts 2^1 = 2\ts 2$ matrix, $\rho_1,$ thus 
\bes
\ \rho_1
\ =
\bbm
 \rho_{11} &  \rho_{12} \\
\rho_{21}  &  \rho_{22}
\ebm
\ees

Note that $\rho_1$ is hermitian and the sum of its diagonal elements, $tr(\rho_1) = 1.$ Using the above mentioned 
standard result, namely, ``Every  $2\ts 2$ hermitian matrix can be uniquely expressed 
as linear combination of Pauli matrices over reals, $\mathbb R$'' it is easy to check that $\rho_1$ can be expressed 
as $$\rho_1 = (\frac{I}{2}+\alpha_1\sigma_1+\alpha_2\sigma_2+\alpha_3\sigma_3)$$ $\alpha_i \in \mathbb{R}, i\in \{1,2,3\}.$ 
It is easy to check that if a matrix expressed as above is hermitian, its trace is equal to one, and further if it is 
positive semidefinite then it will correspond to density matrix of 1-qubit states. Using this representation for density 
matrix of a 1-qubit state we will make the following definitions for ``separable pure'' and ``separable mixed'' states and 
develop the criteria to test the entanglement status in terms of a systematic procedure. Note that all the $\alpha^s$ 
appearing in all the definitions given below are some real constants. 

{\bf Definition 3.3:} A 2-qubit ``pure'' state, represented by the density matrix $\rho_2$ is separable if and only if it can 
be represented as $$\rho_2 = \rho_1^1\ots \rho_1^2 = (\frac{I}{2}+\alpha_1^1\sigma_1+\alpha_2^1\sigma_2+\alpha_3^1\sigma_3)\ots (\frac{I}{2}+\alpha_1^2\sigma_1+\alpha_2^2\sigma_2+\alpha_3^2\sigma_3)$$ 
where $\rho_1^1,\rho_1^2$ are density matrices for 1-qubit states.

{\bf Definition 3.4:} An $N$-qubit ``pure'' state, represented by the density matrix $\rho_N$ is separable if and only if it can 
be represented as $$\rho_N = \prod_{i=1}^{\ots N} \rho_1^i = \prod_{i=1}^{\ots N} (\frac{I}{2}+\alpha_1^i\sigma_1+\alpha_2^i\sigma_2+\alpha_3^i\sigma_3)$$ 
where $\rho_1^i$ are density matrices for 1-qubit states for all $i.$

{\bf Definition 3.5:} A 2-qubit ``mixed'' state, represented by the density matrix $\rho_2$ is separable if and only if it can 
be represented as $$\rho_2 = \sum_{j=1}^M\beta_j(\rho_1^{1j})\ots (\rho_1^{2j}) = \sum_{j=1}^M\beta_j(\frac{I}{2}+\alpha_1^{1j}\sigma_1+\alpha_2^{1j}\sigma_2+\alpha_3^{1j}\sigma_3)\ots (\frac{I}{2}+\alpha_1^{2j}\sigma_1+\alpha_2^{2j}\sigma_2+\alpha_3^{2j}\sigma_3)$$  
for some $M\ge 2$ and $\sum_{j=1}^M\beta_j = 1$ and where tensor products are over density matrices of 1-qubit states. 

{\bf Definition 3.6:} An $N$-qubit ``mixed'' state, represented by the density matrix $\rho_N$ is separable if and only if it can 
be represented as $$\rho_N = \sum_{j=1}^M\beta_j\prod_{i=1}^{\ots N}(\rho_1^{ij}) = \sum_{j=1}^M\beta_j\prod_{i=1}^{\ots N} (\frac{I}{2}+\alpha_1^{ij}\sigma_1+\alpha_2^{ij}\sigma_2+\alpha_3^{ij}\sigma_3)$$ 
for some $M\ge 2$ and $\sum_{j=1}^M\beta_j = 1$ and where tensor products are over density matrices of 1-qubit states. 

\section{A New Criterion and Method for Testing Entanglement status of 2-qubit States}

In this section we will present our new criterion and the method that follows from it for testing entanglement status 
of $N = 2$-qubit states. 

The density operator for $(N = 2)$-qubit states is represented by a $2^2\ts 2^2 = 4\ts 4$ density matrix, $\rho_2$ say, thus, 
\bes
\ \rho_2
\ =
\bbm
\rho_{11}  &  \rho_{12} & \rho_{13} & \rho_{14}  \\
\rho_{21}  &  \rho_{22} & \rho_{23} & \rho_{24}  \\
\rho_{31}  &  \rho_{32} & \rho_{33} & \rho_{34}  \\
\rho_{41}  &  \rho_{42} & \rho_{43} & \rho_{44} 
\ebm
\ees 
We determine whether $\rho_2$ corresponds to a pure state (tr$(\rho_2^2)=1$), or, $\rho_2$ corresponds to a mixed 
state (tr$(\rho_2^2)<1$).

{\bf Case(i)} Suppose $\rho_2$ corresponds to a ``pure'' state. 

We now express $\rho_2$ in terms of Generalized Pauli Basis. Thus, 
$$\rho_2 = \sum_{i,j=0}^3\alpha_{ij}G_{ij}$$ where $\alpha_{ij} = (\rho_2).(G_{ij})/(G_{ij}).(G_{ij})$ and 
$G_{ij} = \sigma_i \ots \sigma_j,$ $\sigma_i,\sigma_j \in \{I (= \sigma_0),\sigma_1,\sigma_2,\sigma_3\}$ 

We form the matrix of coefficients, $\mathbb B = [\alpha_{ij}].$ 

{\bf Theorem 4.1} A 2-qubit ``pure'' quantum state represented by the density matrix $\rho_2$ is separable if and only if 
$\alpha_{11} = \frac{1}{4}$ and the rank of the associated matrix, $\mathbb{B} = [\alpha_{ij}],$ of the coefficients 
that arises when $\rho_2$ is expressed in terms of the Generalized Pauli basis, is equal to one and this unit value of the 
rank makes it possible to express  
$\rho_2 = \rho_1^1\ots \rho_1^2$ where $\rho_1^1\ots \rho_1^2$ is tensor product of density matrices for 1-qubit states.

{\bf Proof:} Let $\rho_2$ be a 2-qubit separable pure state. According to the definition 3.3 
we can express $\rho_2$ in terms of the tensor product of two density matrices for 1-qubit states:  
$$\rho_2 = \rho_1^1\ots \rho_1^2 = (\frac{I}{2}+\alpha_1^1\sigma_1+\alpha_2^1\sigma_2+\alpha_3^1\sigma_3)\ots (\frac{I}{2}+\alpha_1^2\sigma_1+\alpha_2^2\sigma_2+\alpha_3^2\sigma_3)$$ 
 
Carrying out this tensor product and collecting the coefficients of the basis elements of Generalized Pauli basis we can form the matrix $\mathbb{B} = [\alpha_{ij}]$ as
\bes
\ \mathbb{B}
\ =
\bbm
\frac{1}{4}  &  \frac{\alpha_1^2}{2} & \frac{\alpha_2^2}{2} & \frac{\alpha_3^2}{2}\\
\frac{\alpha_1^1}{2}  &  \alpha_1^1\alpha_1^2 & \alpha_1^1\alpha_2^2 & \alpha_1^1\alpha_3^2  \\
\frac{\alpha_2^1}{2}  &  \alpha_2^1\alpha_1^2 & \alpha_2^1\alpha_2^2 & \alpha_2^1\alpha_3^2  \\
\frac{\alpha_3^1}{2}  &  \alpha_3^1\alpha_1^2 & \alpha_3^1\alpha_2^2 & \alpha_3^1\alpha_3^2   
\ebm
\ees 
Clearly, $\alpha_{11} = \frac{1}{4}$ and it is easy to check that the rank of $\mathbb{B}$ is equal to one as required.

Conversely, let $\rho_2$ be a 2-qubit pure state and suppose the rank of the corresponding matrix formed by the coefficients 
when $\rho_2$ is expressed in terms of the Pauli basis, $\mathbb{B} = [\alpha_{ij}]$ is equal to one. Therefore, the 
rows of $\mathbb{B}$ will be mutually proportional. Let $row_i = \delta_i\ts row_1$ where $i\in \{2,3,4\}.$  Since $\mathbb{B}$ 
corresponds to the density matrix, $\rho_2,$ therefore, $\alpha_{11} = \frac{1}{4}.$ It is now easy to check 
that we can express the density matrix, $\rho_2,$ as follows: $$\rho_2 = \rho_1^1\ots \rho_1^2 = (\frac{I}{2}+\frac{\delta_2}{2}\sigma_1+\frac{\delta_3}{2}\sigma_2+\frac{\delta_4}{2}\sigma_3)\ots (\frac{I}{2}+2\alpha_{12}\sigma_1+2\alpha_{13}\sigma_2+2\alpha_{14}\sigma_3)$$ 
where $\alpha_{1j}$ are the matrix elements of the first row of the matrix $\mathbb{B}$ and $\delta^s \in \mathbb{R}.$ 
Therefore, as per definition 3.3 given above this pure state represented by the density matrix, $\rho_2,$ is separable.

{\bf Case(ii)} Suppose $\rho_2$ corresponds to a ``mixed'' state. 

We now express $\rho_2$ in terms of Generalized Pauli Basis. Thus, 
$$\rho_2 = \sum_{i,j=0}^3\alpha_{ij}G_{ij}$$ where $\alpha_{ij} = (\rho_2).(G_{ij})/(G_{ij}).(G_{ij})$ and 
$G_{ij} = \sigma_i \ots \sigma_j,$ $\sigma_i,\sigma_j \in \{I (= \sigma_0),\sigma_1,\sigma_2,\sigma_3\}$ 

As stated in section 1, the matrix $\mathbb B = [\alpha_{ij}].$

{\bf Theorem 4.2} A 2-qubit mixed quantum state represented by the density matrix $\rho_2$ is separable if and only if 
the associated matrix of coefficients that arises when $\rho_2$ is expressed in the 
Pauli basis, $\mathbb{B} = [\alpha_{ij}],$ can be split into sum of matrices $\sum_{k=1}^M\beta_k\mathbb{B}_k,$ 
where $\sum_{k=1}^M\beta_k = 1,$ $\alpha_{11}^k = \frac{1}{4},$ and the rank of each $\mathbb{B}_k = [\alpha_{ij}^k]$ is equal to 
one for all $k.$ 

In other words, a 2-qubit mixed quantum state represented by the density matrix $\rho_2$ is separable if and only if 
the corresponding matrix of coefficients, $\mathbb{B} = [\alpha_{ij}],$ that arises when $\rho_2$ is expressed in the 
Pauli basis can be split into sum of matrices as follows: 
$$\mathbb{B} = [\alpha_{ij}] = \sum_{k=1}^M\beta_k\mathbb{B}_k = \sum_{k=1}^M\beta_k[\alpha_{ij}^k]$$ 
where $\alpha_{11}^k = \frac{1}{4},$ $\sum_{k=1}^M\beta_k = 1,$ and each matrix $\mathbb{B}_k = [\alpha_{ij}^k]$ 
has rank one, for all $k.$ 

{\bf Proof:} Let $\rho_2$ be a 2-qubit separable mixed state. According to definition 3.5 $\rho_2$ can be expressed as sum 
over tensor products of density matrices of 1-qubit states:  
$$\rho_2 = \sum_{k=1}^M\beta_k(\rho_1^{1k}\ots \rho_1^{2k}) = \sum_{k=1}^M\beta_k(\frac{I}{2}+\alpha_1^{1k}\sigma_1+\alpha_2^{1k}\sigma_2+\alpha_3^{1k}\sigma_3)\ots (\frac{I}{2}+\alpha_1^{2k}\sigma_1+\alpha_2^{2k}\sigma_2+\alpha_3^{2k}\sigma_3)$$ 
for some $M\ge 2$ and $\sum_{k=1}^M\beta_k = 1.$ Carrying out the tensor products and collecting the coefficients we can 
form the matrices $\mathbb{B}_k = [\alpha_{ij}^k]$ which are all of the same type as the matrix $\mathbb{B} = [\alpha_{ij}]$ 
given above in the proof of Theorem 4.1 and therefore it is easy to check that rank of each $\mathbb{B}_k$ is equal to 
one as required.

Conversely, let $\rho_2$ be a 2-qubit mixed state and suppose the corresponding matrix of coefficients that arises 
when $\rho_2$ is expressed in the Pauli basis, $\mathbb{B} = [\alpha_{ij}],$ can be split into sum of 
matrices $\sum_{k=1}^M\beta_k\mathbb{B}_k,$ where $\sum_{k=1}^M\beta_k = 1$ and $\alpha_{11}^k = \frac{1}{4}$ and the rank of 
each $\mathbb{B}_k = [\alpha_{ij}^k]$ is equal to one for all $k.$ Therefore, the rows of each $\mathbb{B}_k$ are mutually 
proportional. Therefore, it will be possible to associate a tensor product of type 
$$(\frac{I}{2}+\frac{\delta_2^k}{2}\sigma_1+\frac{\delta_3^k}{2}\sigma_2+\frac{\delta_4^k}{2}\sigma_3)\ots (\frac{I}{2}+2\alpha_{12}^k\sigma_1+2\alpha_{13}^k\sigma_2+2\alpha_{14}^k\sigma_3)$$ 
with each $\mathbb{B}_k$ and using the equation 
$$\mathbb{B} = [\alpha_{ij}] = \sum_{k=1}^M\beta_k\mathbb{B}_k = \sum_{k=1}^M\beta_k[\alpha_{ij}^k]$$ 
we will get 
$$\rho_2 = \sum_{k=1}^M\beta_k\rho_2^k = \sum_{k=1}^M\beta_k(\frac{I}{2}+\frac{\delta_2^k}{2}\sigma_1+\frac{\delta_3^k}{2}\sigma_2+\frac{\delta_4^k}{2}\sigma_3)\ots (\frac{I}{2}+2\alpha_{12}^k\sigma_1+2\alpha_{13}^k\sigma_2+2\alpha_{14}^k\sigma_3)$$ 
implying that $\rho_2$ is separable state.

{\bf Note:} As stated in the introduction we have seen that for $(N=2)$-qubit case we get $4^{(N-2)} = 4^0 =1$ 
matrix, $\mathbb{B} = [\alpha_{ij}],$ called {\bf basic matrix}, made up of coefficients of the basis vectors, $G_{ij,}$ 
that arise while expressing the density matrix, $\rho_2,$ as the unique linear combination in terms of the Generalized Pauli basis. 
We will now proceed to see that for general $N$-qubit case along with certain {\bf basic matrix} of size $4\ts 4$, 
which is similar to $\mathbb{B} = [\alpha_{ij}],$ there arise other matrices. In total there arise $4^{(N-2)}$ matrices 
of size $4\ts 4.$ To understand 
how these matrices arise we think that it is better to consider in details the case 
for $(N=3)$-qubits and to know and record those $4^{(N-2)} = 4^{(3-2)} = 4^1 = 4$ matrices. So we now proceed with 

\section{A New Criterion and Method for Testing Entanglement Status of 3-qubit States}

The density operator for $(N = 3)$-qubit states is represented by a $2^3\ts 2^3 = 8\ts 8$ density matrix, $\rho_3.$ The 
Generalized Pauli basis, for $(N = 3)$-qubit case contains in all 64 basis elements, $G_{ijk} = \sigma_i\ots \sigma_j\ots \sigma_k,$ 
where $\sigma_i, \sigma_j, \sigma_k \in \{I(=\sigma_0), \sigma_1,\sigma_2,\sigma_3\}.$ All basis elements $G_{ijk},$ are 
$8\ts 8$ matrices. 

We first determine whether $\rho_3$ under consideration corresponds to a pure state (tr$(\rho_3^2)=1$), or, $\rho_3$ corresponds to a mixed 
state (tr$(\rho_3^2)<1$).

{\bf Case(i)} Suppose $\rho_3$ corresponds to a ``pure'' state. 

We now express $\rho_3$ in terms of Generalized Pauli Basis. Thus, 
$$\rho_3 = \sum_{i,j,k=0}^3\alpha_{ijk}G_{ijk}$$ where 
$G_{ijk} = \sigma_i \ots \sigma_j\ots \sigma_k,$ $\sigma_i,\sigma_j,\sigma_k \in \{I (= \sigma_0),\sigma_1,\sigma_2,\sigma_3\}$ 
and $\alpha_{ijk} = (\rho_3).(G_{ijk})/(G_{ijk}).(G_{ijk}).$

{\bf Theorem 5.1} A 3-qubit pure quantum state represented by the density matrix $\rho_3$ is separable if and only if 
$\alpha_{000} = \frac{1}{8}$ and among the following four $4\ts 4$ matrices of coefficients 
$A_i,i\in \{0,1,2,3\},$ the rank of the matrix $A_0$ (called the {\bf basic matrix}) is equal to one, and $A_j = k_jA_0$ where $k_j$ 
are some constants and $j\in \{1,2,3\}$ (thus, the matrices $A_j,j\in \{1,2,3\}$ are some constant multiples of the basic 
matrix $A_0$ and so obviously their rank is also equal to one). These matrices, $A_i,i\in \{0,1,2,3\}$ are as follows: 
\bes
\ A_0
\ =
\bbm
\alpha_{000}  &  \alpha_{001} & \alpha_{002} & \alpha_{003}  \\
\alpha_{010}  &  \alpha_{011} & \alpha_{012} & \alpha_{013}  \\
\alpha_{020}  &  \alpha_{021} & \alpha_{022} & \alpha_{023}  \\
\alpha_{030}  &  \alpha_{031} & \alpha_{032} & \alpha_{033}  
\ebm
\ees 

\bes
\ A_1
\ =
\bbm
\alpha_{100}  &  \alpha_{101} & \alpha_{102} & \alpha_{103}  \\
\alpha_{110}  &  \alpha_{111} & \alpha_{112} & \alpha_{113}  \\
\alpha_{120}  &  \alpha_{121} & \alpha_{122} & \alpha_{123}  \\
\alpha_{130}  &  \alpha_{131} & \alpha_{132} & \alpha_{133}  
\ebm
\ees

\bes
\ A_2
\ =
\bbm
\alpha_{200}  &  \alpha_{201} & \alpha_{202} & \alpha_{203}  \\
\alpha_{210}  &  \alpha_{211} & \alpha_{212} & \alpha_{213}  \\
\alpha_{220}  &  \alpha_{221} & \alpha_{222} & \alpha_{223}  \\
\alpha_{230}  &  \alpha_{231} & \alpha_{232} & \alpha_{233}  
\ebm
\ees

\bes
\ A_3
\ =
\bbm
\alpha_{300}  &  \alpha_{301} & \alpha_{302} & \alpha_{303}  \\
\alpha_{310}  &  \alpha_{311} & \alpha_{312} & \alpha_{313}  \\
\alpha_{320}  &  \alpha_{321} & \alpha_{322} & \alpha_{323}  \\
\alpha_{330}  &  \alpha_{331} & \alpha_{332} & \alpha_{333}  
\ebm
\ees

{\bf Proof:} Let $\rho_3$ be a 3-qubit separable pure state. According to definition 3.4 $\rho_3$ can be represented as 
$$\rho_3 = \prod_{i=1}^{\ots 3} (\frac{I}{2}+\alpha_1^i\sigma_1+\alpha_2^i\sigma_2+\alpha_3^i\sigma_3)$$ 

Carrying out the tensor product and collecting the coefficients we can form the basic matrix $A_0$ as
\bes
\ A_0
\ =
\bbm
\frac{1}{8}  &  \frac{\alpha_1^3}{4} & \frac{\alpha_2^3}{4} & \frac{\alpha_3^3}{4} \\
\frac{\alpha_1^2}{4}  & \frac{\alpha_1^2\alpha_1^3}{2} & \frac{\alpha_1^2\alpha_2^3}{2} & \frac{\alpha_1^2\alpha_3^3}{2} \\
\frac{\alpha_2^2}{4}  & \frac{\alpha_2^2\alpha_1^3}{2} & \frac{\alpha_2^2\alpha_2^3}{2} & \frac{\alpha_2^2\alpha_3^3}{2} \\
\frac{\alpha_3^2}{4}  & \frac{\alpha_3^2\alpha_1^3}{2} & \frac{\alpha_3^2\alpha_2^3}{2} & \frac{\alpha_3^2\alpha_3^3}{2} \\ 
\ebm
\ees 
Clearly, $\alpha_{000} = \frac{1}{8}$ and it is easy to check that the rank of $A_0$ is equal to one as required. Further, it 
is easy to check that $A_1 = 2\alpha_1^1A_0,$ $A_2 = 2\alpha_2^1A_0,$ $A_3 = 2\alpha_3^1A_0$ as required. 

Conversely, let $\rho_3$ be a 3-qubit pure state and suppose the rank of the associated matrix, $A_0$ is equal to one and 
the matrices $A_j,j\in \{1,2,3\}$ are some constant multiples of the basic 
matrix $A_0$ and so obviously their rank is also equal to one.
Therefore, the rows of $A_0$ are mutually proportional. Let $row_i = \delta_i\ts row_1$ where $i\in \{2,3,4\}.$  
Since $A_0$ corresponds to the density matrix of 3-qubit state $\rho_3,$ therefore, $\alpha_{000} = \frac{1}{8}.$ It is now easy 
to check that we can express the density matrix, $\rho_3,$ as 
follows: $$\rho_3 = (\frac{I}{2}+\frac{\gamma_1}{2}\sigma_1+\frac{\gamma_2}{2}\sigma_2+\frac{\gamma_3}{2}\sigma_3)\ots (\frac{I}{2}+\frac{\delta_2}{2}\sigma_1+\frac{\delta_3}{2}\sigma_2+\frac{\delta_4}{2}\sigma_3)\ots (\frac{I}{2}+4\alpha_{001}\sigma_1+4\alpha_{002}\sigma_2+4\alpha_{003}\sigma_3)$$ 
Therefore, as per definition 3.4 given above this pure state represented by the density matrix, $\rho_3,$ is separable. 

{\bf Note:} $A_0$ is made up of the coefficients of the basis elements in the Generalized Pauli basis that occur in the 
expansion of $$(\frac{I}{2})\ots (\frac{I}{2}+\frac{\delta_2}{2}\sigma_1+\frac{\delta_3}{2}\sigma_2+\frac{\delta_4}{2}\sigma_3)\ots (\frac{I}{2}+4\alpha_{001}\sigma_1+4\alpha_{002}\sigma_2+4\alpha_{003}\sigma_3)$$ 
Similarly, $A_1$ is made up of the coefficients of the basis elements in the Generalized Pauli basis that occur in the 
expansion of $$(\frac{\gamma_1}{2}\sigma_1)\ots (\frac{I}{2}+\frac{\delta_2}{2}\sigma_1+\frac{\delta_3}{2}\sigma_2+\frac{\delta_4}{2}\sigma_3)\ots (\frac{I}{2}+4\alpha_{001}\sigma_1+4\alpha_{002}\sigma_2+4\alpha_{003}\sigma_3)$$ 
$A_2$ is made up of the coefficients of the basis elements in the Generalized Pauli basis that occur in the 
expansion of $$(\frac{\gamma_2}{2}\sigma_2)\ots (\frac{I}{2}+\frac{\delta_2}{2}\sigma_1+\frac{\delta_3}{2}\sigma_2+\frac{\delta_4}{2}\sigma_3)\ots (\frac{I}{2}+4\alpha_{001}\sigma_1+4\alpha_{002}\sigma_2+4\alpha_{003}\sigma_3)$$ 
$A_3$ is made up of the coefficients of the basis elements in the Generalized Pauli basis that occur in the 
expansion of $$(\frac{\gamma_3}{2}\sigma_3)\ots (\frac{I}{2}+\frac{\delta_2}{2}\sigma_1+\frac{\delta_3}{2}\sigma_2+\frac{\delta_4}{2}\sigma_3)\ots (\frac{I}{2}+4\alpha_{001}\sigma_1+4\alpha_{002}\sigma_2+4\alpha_{003}\sigma_3)$$ 

{\bf Case(ii)} Suppose $\rho_3$ corresponds to a ``mixed'' state. 

We now express $\rho_3$ in terms of Generalized Pauli Basis. Thus, 
$$\rho_3 = \sum_{i,j,k=0}^3\alpha_{ijk}G_{ijk}$$ where $G_{ijk} = \sigma_i \ots \sigma_j\ots \sigma_k,$ and $\alpha_{ijk} = (\rho_3).(G_{ijk})/(G_{ijk}).(G_{ijk})$ 
$\sigma_i,\sigma_j,\sigma_k \in \{I (= \sigma_0),\sigma_1,\sigma_2,\sigma_3\}$  

{\bf Theorem 5.2} A 3-qubit mixed quantum state represented by the density matrix $\rho_3$ is separable if and only if 
the associated basic matrix, $A_0,$ of coefficients that arises when $\rho_3$ is expressed in the 
Pauli basis can be split into sum of matrices $\sum_{k=1}^M\beta_kA_0^k,$ 
where $\sum_{k=1}^M\beta_k = 1$ and $\alpha_{000}^k = \frac{1}{8}$ and the rank of each $A_0^k$ is equal to 
one for all $k$ and the other matrices $A_j,j\in \{1,2,3\}$ are some constant multiples of the basic 
matrix $A_0.$ 

In other words, a 3-qubit mixed quantum state represented by the density matrix $\rho_3$ is separable if and only if 
the associated basic matrix of coefficients, $A_0,$ that arises when $\rho_3$ is expressed in the 
Pauli basis can be split into sum of matrices as follows: 
$$A_0 = \sum_{k=1}^M\beta_k A_0^k$$ 
where each matrix $A_0^k$ has rank one and $\alpha_{000}^k = \frac{1}{8}$ for all $k.$ Further, $\sum_{k=1}^M\beta_k = 1$ 
and the other matrices $A_j,j\in \{1,2,3\}$ are some constant multiples of the basic matrix $A_0.$ 

{\bf Proof:} The proof will be similar to the proofs of the Theorem 4.2 and Theorem 5.1 given above and can be easily 
obtained by proceeding on similar lines.

\section{Algorithms:} When a quantum state containing $N$-qubits is given in terms of its corresponding density matrix, $\rho_N,$ 
it is important to know for various applications whether this state is ``entangled'' or ``separable''. 

It turns out, whatever the dimensions of the density operator, $\rho_N,$ that $(\rho_N)^2 = \rho_N$ and 
$tr((\rho_N)^2) = 1$ if and only if $\rho_N$ is the density operator corresponding to a pure state. If the state described 
by $\rho_N$ is not pure, but is instead mixed, then $(\rho_N)^2 \neq \rho_N$ and $tr((\rho_N)^2) < 1.$ These properties 
can be used to decide whether a given state is pure or mixed.

Before proceeding to check its entanglement status it is useful to know whether the state defined in terms of the density 
matrix $\rho_N$ corresponds to a ``pure'' state or a ``mixed'' state and as stated above there is available an easy test 
for it. One simply needs to determine the trace, $tr(\rho_N^2).$ 

(i) The state is ``pure'' if and only if $tr(\rho_N^2) = 1.$ 

(ii) The state is ``mixed'' if and only if $tr(\rho_N^2) < 1.$ 

We now proceed to develop algorithms to test the entanglement status of 2-qubit states given in terms of their corresponding 
density matrix, $\rho_2.$

{\bf Algorithms for 2-qubit states:}

We first determine using the above given criteria whether the state under consideration is ``pure'' or ``mixed''.

{\bf Case(i)} Suppose $\rho_2$ corresponds to a ``pure'' state.

{\bf Algorithm 1:}

{\bf Step 1:} Express $\rho_2$ in terms of Generalized Pauli Basis. Thus, 
$$\rho_2 = \sum_{i,j=0}^3\alpha_{ij}G_{ij}$$ where $\alpha_{ij} = (\rho_2).(G_{ij})/(G_{ij}).(G_{ij})$ and 
$G_{ij} = \sigma_i \ots \sigma_j,$ $\sigma_i,\sigma_j \in \{I (= \sigma_0),\sigma_1,\sigma_2,\sigma_3\}$ 
and form the matrix $\mathbb B = [\alpha_{ij}].$ 

{\bf Step 2:} Check whether the rank of $\mathbb B$ is equal to one or not. 

(a) If the rank of $\mathbb B$ is equal to one then declare that the state is separable and find out the expression 
for $\rho_2$ as $\rho_2 = \rho_1^A\ots \rho_1^B$ using Theorem 4.1.

(b) If the rank of $\mathbb B$ is not equal to one then declare that the state is entangled. 

{\bf Case(ii)} Suppose $\rho_2$ corresponds to a ``mixed`` state. 

{\bf Algorithm 2:}

{\bf Step 1:} Express $\rho_3$ in terms of Generalized Pauli Basis. Thus, 
$$\rho_2 = \sum_{i,j=0}^3\alpha_{ij}G_{ij}$$ where $\alpha_{ij} = (\rho_2).(G_{ij})/(G_{ij}).(G_{ij})$ and 
$G_{ij} = \sigma_i \ots \sigma_j,$ $\sigma_i,\sigma_j \in \{I (= \sigma_0),\sigma_1,\sigma_2,\sigma_3\}$ 
and form the matrix $\mathbb B = [\alpha_{ij}].$ 

{\bf Step 2:} Check whether the rank of the sub-matrix of $\mathbb B$ formed by the second, third, and fourth rows is 
equal to one and whether it is possible to split the elements of the first row of $\mathbb B$ as $\alpha_{1j} = \alpha_{1j}^1 + \alpha_{1j}^2$ and rewrite the 
matrix, $\mathbb B,$ as $\mathbb B = \mathbb B^1 + \mathbb B^2$ where
\bes
\ \mathbb{B}^1
\ =
\bbm
\alpha_{11}^1  &  \alpha_{12}^1 & \alpha_{13}^1 & \alpha_{14}^1  \\
\alpha_{21}  &  \alpha_{22} & \alpha_{23} & \alpha_{24}  \\
\alpha_{31}  &  \alpha_{32} & \alpha_{33} & \alpha_{34}  \\
\alpha_{41}  &  \alpha_{42} & \alpha_{43} & \alpha_{44} 
\ebm
\ees 
and 
\bes
\ \mathbb{B}^2
\ =
\bbm
\alpha_{11}^2  &  \alpha_{12}^2 & \alpha_{13}^2 & \alpha_{14}^2  \\
      0        &        0       &       0       &       0        \\ 
      0        &        0       &       0       &       0        \\
      0        &        0       &       0       &       0       
\ebm
\ees   
such that now the rank of the matrix $\mathbb{B}^1$ is indeed equal to one (note that the rank of the matrix $\mathbb{B}^2$ 
is obviously equal to one).

(c) If yes, then declare using Theorem 4.2 that the mixed state $\rho_2$ is separable and then one can easily express $\rho_2$ as 
$\rho_2 = \beta_1(\rho_1^A\ots \rho_1^B) + \beta_2(\rho_1^C\ots \rho_1^D),$ where $\beta_1 + \beta_2 = 1.$ 

(d) If no, then proceed to 

{\bf Step 3:} Check whether the rank of the sub-matrix of $\mathbb B$ formed by any two rows, say the second and third rows 
is equal to one and whether it is possible to split the elements of the first row of $\mathbb B$ 
as $\alpha_{1j} = \alpha_{1j}^1 + \alpha_{1j}^2 + \alpha_{1j}^3$ and rewrite the 
matrix, $\mathbb B,$ as $\mathbb B = \mathbb B^1 + \mathbb B^2 + \mathbb B^3$ where
\bes
\ \mathbb{B}^1
\ =
\bbm
\alpha_{11}^1  &  \alpha_{12}^1 & \alpha_{13}^1 & \alpha_{14}^1  \\
\alpha_{21}  &  \alpha_{22} & \alpha_{23} & \alpha_{24}  \\
\alpha_{31}  &  \alpha_{32} & \alpha_{33} & \alpha_{34}  \\
      0      &        0     &       0     &       0 
\ebm
\ees 

\bes
\ \mathbb{B}^2
\ =
\bbm
\alpha_{11}^2  &  \alpha_{12}^2 & \alpha_{13}^2 & \alpha_{14}^2  \\
      0        &        0       &       0       &       0        \\ 
      0        &        0       &       0       &       0        \\
\alpha_{41}  &  \alpha_{42} & \alpha_{43} & \alpha_{44}     
\ebm
\ees
and
\bes
\ \mathbb{B}^3
\ =
\bbm
\alpha_{11}^3  &  \alpha_{12}^3 & \alpha_{13}^3 & \alpha_{14}^3  \\
      0        &        0       &       0       &       0        \\ 
      0        &        0       &       0       &       0        \\
      0        &        0       &       0       &       0       
\ebm
\ees  
such that now the rank of the matrix $\mathbb{B}^1,\mathbb{B}^2$ is indeed equal to one 
(note that the rank of the matrix $\mathbb{B}^3$ is obviously equal to one).

(e) If yes, then declare using Theorem 4.2 that the mixed state $\rho_2$ is separable and then one can easily express $\rho_2$ as 
$\rho_2 = \beta_1(\rho_1^A\ots \rho_1^B) + \beta_2(\rho_1^C\ots \rho_1^D + \beta_3(\rho_1^E\ots \rho_1^F),$ 
where $\beta_1 + \beta_2 + \beta_3 = 1.$ 

(f) If no, then proceed to 

{\bf Step 4:} Suppose now that the rank of the sub-matrix of $\mathbb B$ formed by any two of the two rows is not 
equal to one. Check now whether it is possible to split the elements of the first row of $\mathbb B$ 
as $\alpha_{1j} = \alpha_{1j}^1 + \alpha_{1j}^2 + \alpha_{1j}^3 + \alpha_{1j}^4$ and we can rewrite the 
matrix, $\mathbb B,$ as $\mathbb B = \mathbb B^1 + \mathbb B^2 + \mathbb B^3 + \mathbb B^4$ where
\bes
\ \mathbb{B}^1
\ =
\bbm
\alpha_{11}^1  &  \alpha_{12}^1 & \alpha_{13}^1 & \alpha_{14}^1  \\
\alpha_{21}  &  \alpha_{22} & \alpha_{23} & \alpha_{24}  \\
      0      &        0     &       0     &       0      \\
      0      &        0     &       0     &       0 
\ebm
\ees

\bes
\ \mathbb{B}^2
\ =
\bbm
\alpha_{11}^2  &  \alpha_{12}^2 & \alpha_{13}^2 & \alpha_{14}^2  \\
      0      &        0     &       0     &       0      \\
\alpha_{31}  &  \alpha_{32} & \alpha_{33} & \alpha_{34}  \\
      0      &        0     &       0     &       0 
\ebm
\ees 

\bes
\ \mathbb{B}^3
\ =
\bbm
\alpha_{11}^3  &  \alpha_{12}^3 & \alpha_{13}^3 & \alpha_{14}^3  \\
      0        &        0       &       0       &       0        \\ 
      0        &        0       &       0       &       0        \\
\alpha_{41}  &  \alpha_{42} & \alpha_{43} & \alpha_{44}     
\ebm
\ees
and
\bes
\ \mathbb{B}^4
\ =
\bbm
\alpha_{11}^4  &  \alpha_{12}^4 & \alpha_{13}^4 & \alpha_{14}^4  \\
      0        &        0       &       0       &       0        \\ 
      0        &        0       &       0       &       0        \\
      0        &        0       &       0       &       0       
\ebm
\ees  
such that now the rank of the matrix $\mathbb{B}^1,\mathbb{B}^2,\mathbb{B}^3$ is indeed equal to one 
(note that the rank of the matrix $\mathbb{B}^4$ is obviously equal to one). 

(g) If yes, then declare using Theorem 4.2 that the mixed state $\rho_2$ is separable and then one can easily express $\rho_2$ as 
$\rho_2 = \beta_1(\rho_1^A\ots \rho_1^B) + \beta_2(\rho_1^C\ots \rho_1^D + \beta_3(\rho_1^E\ots \rho_1^F) + \beta_4(\rho_1^G\ots \rho_1^H),$ 
where $\beta_1 + \beta_2 + \beta_3 + \beta_4 = 1.$ 

(h) If no, then as per Theorem 4.2 the state $\rho_2$ is not separable and so declare that the state is entangled.  

{\bf Algorithms for 3-qubit states:}

We first determine using the above given criteria whether the state under consideration is ``pure'' or ``mixed''.

{\bf Case(i)} Suppose $\rho_3$ corresponds to a ``pure'' state.

{\bf Algorithm 3:}

{\bf Step 1:} Express $\rho_3$ in terms of Generalized Pauli Basis. Thus, 
$$\rho_3 = \sum_{i,j,k=0}^3\alpha_{ijk}G_{ijk}$$ where $\alpha_{ijk} = (\rho_3).(G_{ijk})/(G_{ijk}).(G_{ijk})$ and 
$G_{ijk} = \sigma_i \ots \sigma_j\ots \sigma_k,$ $\sigma_i,\sigma_j,\sigma_k \in \{I (= \sigma_0),\sigma_1,\sigma_2,\sigma_3\}$ 
and form the matrices $A_0,A_1,A_2,A_3$ 

{\bf Step 2:} Check whether the rank of $A_0$ is equal to one or not. 

(a) If the rank of $A_0$ is equal to one and $A_i = k_i\ts A_0, i \in \{(1,2,3)\}$ where $k_i$ are some constants, then 
declare that the state is separable and find out the expression 
for $\rho_3$ as $\rho_3 = \rho_1^A\ots \rho_1^B\ots \rho_1^C$ using Theorem 5.1.

(b) If the rank of $A_0$ is not equal to one then declare that the state is entangled. 

(c) If the rank of $A_0$ is equal to one but $A_i \neq k_i\ts A_0$ for some one or more $i \in \{(1,2,3)\}$ where $k_i$ are some constants, then 
declare that the state is entangled.  

{\bf Case(ii)} Suppose $\rho_3$ corresponds to a ``mixed`` state. 

{\bf Algorithm 4:}

{\bf Step 1:} Express $\rho_3$ in terms of Generalized Pauli Basis. Thus, 
$$\rho_3 = \sum_{i,j,k=0}^3\alpha_{ijk}G_{ijk}$$ where $\alpha_{ijk} = (\rho_3).(G_{ijk})/(G_{ijk}).(G_{ijk})$ and 
$G_{ijk} = \sigma_i \ots \sigma_j\ots \sigma_k,$ $\sigma_i,\sigma_j,\sigma_k \in \{I (= \sigma_0),\sigma_1,\sigma_2,\sigma_3\}$ 
and form the matrices $A_0,A_1,A_2,A_3$ and check whether $A_i = k_i\ts A_0, i \in \{(1,2,3)\}$ where $k_i$ are some constants. 

{\bf Case (1):} Suppose $A_i = k_i\ts A_0, i \in \{(1,2,3)\}$ where $k_i$ are some constants.

{\bf Step 2:} Check whether the rank of the sub-matrix of $A_0$ formed by the second, third, and fourth rows is 
equal to one and whether it is possible to split the elements of the first row of $A_0$ as $\alpha_{00j} = \alpha_{00j}^1 + \alpha_{00j}^2$ and rewrite the 
matrix, $A_0$ as $A_0 = A_0^1 + A_0^2$ where
\bes
\ A_0^1
\ =
\bbm
\alpha_{000}^1  &  \alpha_{001}^1 & \alpha_{002}^1 & \alpha_{003}^1  \\
\alpha_{010}  &  \alpha_{011} & \alpha_{012} & \alpha_{013}  \\
\alpha_{020}  &  \alpha_{021} & \alpha_{022} & \alpha_{023}  \\
\alpha_{030}  &  \alpha_{031} & \alpha_{032} & \alpha_{033} 
\ebm
\ees 
and 
\bes
\ A_0^2
\ =
\bbm
\alpha_{000}^2  &  \alpha_{001}^2 & \alpha_{002}^2 & \alpha_{003}^2  \\
      0        &        0       &       0       &       0        \\ 
      0        &        0       &       0       &       0        \\
      0        &        0       &       0       &       0       
\ebm
\ees   
such that now the rank of the matrix $A_0^1$ is indeed equal to one (note that the rank of the matrix $A_0^2$ 
is obviously equal to one).

(c) If yes, then declare using Theorem 5.2 that the mixed state $\rho_3$ is separable and then one can easily express $\rho_3$ as 
$\rho_3 = \beta_1(\rho_1^A\ots \rho_1^B\ots \rho_1^C) + \beta_2(\rho_1^D\ots \rho_1^E\ots \rho_1^F),$ where $\beta_1 + \beta_2 = 1.$ 

(d) If no, then proceed to 

{\bf Step 3:} Check whether the rank of the sub-matrix of $A_0$ formed by any two rows, say the second and third rows 
is equal to one and whether it is possible to split the elements of the first row of $A_0$ 
as $\alpha_{00j} = \alpha_{00j}^1 + \alpha_{00j}^2 + \alpha_{00j}^3$ and rewrite the 
matrix, $A_0$ as $A_0 = A_0^1 + A_0^2 + A_0^3$ where
\bes
\ A_0^1
\ =
\bbm
\alpha_{000}^1  &  \alpha_{001}^1 & \alpha_{002}^1 & \alpha_{003}^1  \\
\alpha_{010}  &  \alpha_{011} & \alpha_{012} & \alpha_{013}  \\
\alpha_{020}  &  \alpha_{021} & \alpha_{022} & \alpha_{023}  \\
      0      &        0     &       0     &       0 
\ebm
\ees 

\bes
\ A_0^2
\ =
\bbm
\alpha_{000}^2  &  \alpha_{001}^2 & \alpha_{002}^2 & \alpha_{003}^2  \\
      0        &        0       &       0       &       0        \\ 
      0        &        0       &       0       &       0        \\
\alpha_{030}  &  \alpha_{031} & \alpha_{032} & \alpha_{033}     
\ebm
\ees
and
\bes
\ A_0^3
\ =
\bbm
\alpha_{000}^3  &  \alpha_{001}^3 & \alpha_{002}^3 & \alpha_{003}^3  \\
      0        &        0       &       0       &       0        \\ 
      0        &        0       &       0       &       0        \\
      0        &        0       &       0       &       0       
\ebm
\ees  
such that now the rank of the matrix $A_0^1,A_0^2$ is indeed equal to one 
(note that the rank of the matrix $A_0^3$ is obviously equal to one).

(e) If yes, then declare using Theorem 5.2 that the mixed state $\rho_3$ is separable and then one can easily express $\rho_3$ as 
$\rho_3 = \beta_1(\rho_1^A\ots \rho_1^B\ots \rho_1^C) + \beta_2(\rho_1^D\ots \rho_1^E\ots \rho_1^F) + \beta_3(\rho_1^G\ots \rho_1^H\ots \rho_1^I),$ 
where $\beta_1 + \beta_2 + \beta_3 = 1.$ 

(f) If no, then proceed to 

{\bf Step 4:} Suppose now that the rank of the sub-matrix of $A_0$ formed by any two rows out of second,third and fourth rows is not 
equal to one. Check now whether it is possible to split the elements of the first row of $A_0$ 
as $\alpha_{00j} = \alpha_{00j}^1 + \alpha_{00j}^2 + \alpha_{00j}^3 + \alpha_{00j}^4$ and we can rewrite the 
matrix, $A_0$ as $A_0 = A_0^1 + A_0^2 + A_0^3 + A_0^4$ where
\bes
\ A_0^1
\ =
\bbm
\alpha_{000}^1  &  \alpha_{001}^1 & \alpha_{002}^1 & \alpha_{003}^1  \\
\alpha_{010}  &  \alpha_{011} & \alpha_{012} & \alpha_{013}  \\
      0      &        0     &       0     &       0      \\
      0      &        0     &       0     &       0 
\ebm
\ees

\bes
\ A_0^2
\ =
\bbm
\alpha_{000}^2  &  \alpha_{001}^2 & \alpha_{002}^2 & \alpha_{003}^2  \\
      0      &        0     &       0     &       0      \\
\alpha_{020}  &  \alpha_{021} & \alpha_{022} & \alpha_{023}  \\
      0      &        0     &       0     &       0 
\ebm
\ees 

\bes
\ A_0^3
\ =
\bbm
\alpha_{000}^3  &  \alpha_{001}^3 & \alpha_{002}^3 & \alpha_{003}^3  \\
      0        &        0       &       0       &       0        \\ 
      0        &        0       &       0       &       0        \\
\alpha_{030}  &  \alpha_{031} & \alpha_{032} & \alpha_{033}     
\ebm
\ees
and
\bes
\ A_0^4
\ =
\bbm
\alpha_{000}^4  &  \alpha_{001}^4 & \alpha_{002}^4 & \alpha_{003}^4  \\
      0        &        0       &       0       &       0        \\ 
      0        &        0       &       0       &       0        \\
      0        &        0       &       0       &       0       
\ebm
\ees  
such that now the rank of the matrix $A_0^1,A_0^2,A_0^3$ is indeed equal to one 
(note that the rank of the matrix $A_0^4$ is obviously equal to one). 

(g) If yes, then declare using Theorem 5.2 that the mixed state $\rho_3$ is separable and then one can easily express $\rho_3$ as 
$\rho_3 = \beta_1(\rho_1^A\ots \rho_1^B\ots \rho_1^C) + \beta_2(\rho_1^D\ots \rho_1^E\ots \rho_1^F) + \beta_3(\rho_1^G\ots \rho_1^H\ots \rho_1^I) + \beta_4(\rho_1^J\ots \rho_1^K\ots \rho_1^L),$ 
where $\beta_1 + \beta_2 + \beta_3 + \beta_4 = 1.$ 

(h) If no, then as per Theorem 5.2 the state $\rho_3$ is not separable and so declare that the state is entangled.  

{\bf Case (2):} Suppose $A_i \neq k_i\ts A_0$ for some one or more $i \in \{(1,2,3)\}$ where $k_i$ are some constants, then 
declare that the state is entangled.  

\section{Further Generalizations:} It is possible for one to continue on these lines and to deal with the higher dimensional 
cases by proceeding on similar lines. One will need to form the associated $4\ts 4$ matrices of coefficients and check 
whether they posses required properties  as seen in the above cases. As mentioned above the number of these matrices grow with 
the dimensions as they are $4^{(N-2)}$ in number for the $N$ dimensional case. 

Are the matrices that arise in $N$ dimensional case related to the matrices that arise in $(N-1)$ dimensional case? Yes.  

We give below this relation for the case of $N = 4$ with the case of $N = 3,$ i.e. we will relate these 16 matrices, 
say $B_0,B_1,\cdots,B_{15},$ that arise for the case $N = 4$ with the four matrices $A_0,A_1,A_2,A_3$ that aroused for 
the case of $N = 3.$ 

(1) The matrix $B_0$ is the one obtained by attaching additional suffix ``0'' from the left side to each element of the 
above matrix $A_0,$ e.g. an element, say $\alpha_{0bc} \to \alpha_{00bc}.$ Hence etc.

(2) The matrix $B_1$ is the one obtained by attaching additional suffix ``0'' from the left side to each element of the 
above matrix $A_1,$ e.g. an element, say $\alpha_{1bc} \to \alpha_{01bc}.$ Hence etc.

(3) The matrix $B_2$ is the one obtained by attaching additional suffix ``0'' from the left side to each element of the 
above matrix $A_2,$ e.g. an element, say $\alpha_{2bc} \to \alpha_{02bc}.$ Hence etc.

(4) The matrix $B_3$ is the one obtained by attaching additional suffix ``0'' from the left side to each element of the 
above matrix $A_3,$ e.g. an element, say $\alpha_{3bc} \to \alpha_{03bc}.$ Hence etc.

(5) The matrix $B_4$ is the one obtained by attaching additional suffix ``1'' from the left side to each element of the 
above matrix $A_0,$ e.g. an element, say $\alpha_{0bc} \to \alpha_{10bc}.$ Hence etc.

(6) The matrix $B_5$ is the one obtained by attaching additional suffix ``1'' from the left side to each element of the 
above matrix $A_1,$ e.g. an element, say $\alpha_{1bc} \to \alpha_{11bc}.$ Hence etc.

(7) The matrix $B_6$ is the one obtained by attaching additional suffix ``1'' from the left side to each element of the 
above matrix $A_2,$ e.g. an element, say $\alpha_{2bc} \to \alpha_{12bc}.$ Hence etc.

(8) The matrix $B_7$ is the one obtained by attaching additional suffix ``1'' from the left side to each element of the 
above matrix $A_3,$ e.g. an element, say $\alpha_{3bc} \to \alpha_{13bc}.$ Hence etc.

(9) The matrix $B_8$ is the one obtained by attaching additional suffix ``2'' from the left side to each element of the 
above matrix $A_0,$ e.g. an element, say $\alpha_{0bc} \to \alpha_{20bc}.$ Hence etc.

(10) The matrix $B_9$ is the one obtained by attaching additional suffix ``2'' from the left side to each element of the 
above matrix $A_1,$ e.g. an element, say $\alpha_{1bc} \to \alpha_{21bc}.$ Hence etc.

(11) The matrix $B_{10}$ is the one obtained by attaching additional suffix ``2'' from left side to each element of the 
above matrix $A_2,$ e.g. an element, say $\alpha_{2bc} \to \alpha_{22bc}.$ Hence etc.

(12) The matrix $B_{11}$ is the one obtained by attaching additional suffix ``2'' from the left side to each element of 
the above matrix $A_3,$ e.g. an element, say $\alpha_{3bc} \to \alpha_{23bc}.$ Hence etc.

(13) The matrix $B_{12}$ is the one obtained by attaching additional suffix ``3'' from the left side to each element of 
the above matrix $A_0,$ e.g. an element, say $\alpha_{0bc} \to \alpha_{30bc}.$ Hence etc.

(14) The matrix $B_{13}$ is the one obtained by attaching additional suffix ``3'' from the left side to each element of 
the above matrix $A_1,$ e.g. an element, say $\alpha_{1bc} \to \alpha_{31bc}.$ Hence etc.

(15) The matrix $B_{14}$ is the one obtained by attaching additional suffix ``3'' from the left side to each element of 
the above matrix $A_2,$ e.g. an element, say $\alpha_{2bc} \to \alpha_{32bc}.$ Hence etc.

(16) The matrix $B_{15}$ is the one obtained by attaching additional suffix ``3'' from the left side to each element of 
the above matrix $A_3,$ e.g. an element, say $\alpha_{3bc} \to \alpha_{33bc}.$ Hence etc.

\section{Examples:} 

{\bf Example 1:} We consider a 2-qubit quantum state defined by the following density matrix, $\rho_2,$ as given below:
\bes
\ \rho_2
\ =
\bbm
\frac{1}{36}  &  \frac{1}{9\sqrt{2}} & \frac{1}{12\sqrt{3}} & \frac{1}{3\sqrt{6}}  \\
\frac{1}{9\sqrt{2}}  &  \frac{2}{9} & \frac{1}{3\sqrt{6}} & \frac{2}{3\sqrt{3}}  \\
\frac{1}{12\sqrt{3}}  &  \frac{1}{3\sqrt{6}} & \frac{1}{12} & \frac{1}{3\sqrt{2}}  \\
\frac{1}{3\sqrt{6}}  &  \frac{2}{3\sqrt{3}} & \frac{1}{3\sqrt{2}} & \frac{2}{3} 
\ebm
\ees 
It is easy to see that $tr((\rho_2)^2) = 1$ therefore $\rho_2$ represents a ``pure'' state.  

We now express $\rho_2$ in terms of Generalized Pauli basis and construct the matrix of 
coefficients, $\mathbb{B},$ as given below: 
\bes
\ \mathbb{B}
\ =
\bbm
\frac{1}{4}  &  \frac{\sqrt{2}}{9} & 0 & -\frac{7}{36}  \\
\frac{\sqrt{3}}{8}  &  \frac{\sqrt{6}}{18} & 0 & -\frac{7\sqrt{3}}{72}  \\
0  &  0 & 0 & 0  \\
-\frac{1}{8}  &  -\frac{\sqrt{2}}{18} & 0 & \frac{7}{72} 
\ebm
\ees 
and apply Theorem 4.1 to it. Note that $\alpha_{11} = \frac{1}{4}.$
It is easy to check that the rank of $\mathbb{B}$ is equal to one therefore the 2-qubit quantum state defined by the 
density matrix, $\rho_2,$ is a {\bf pure separable} state. We can now  
easily check that (i) $row_2 = \frac{\sqrt{3}}{2}\ts row_1,$ (ii) $row_3 = 0\ts row_1,$ 
and (iii) $row_4 = -\frac{1}{2}\ts row_1.$ Therefore, following Algorithm 1 we have 
$$\rho_2 = (\frac{I}{2}+\frac{\sqrt{3}}{4}\sigma_1-\frac{1}{4}\sigma_3)\ots (\frac{I}{2}+\frac{2\sqrt{2}}{9}\sigma_1-\frac{7}{18}\sigma_3)$$

It is easy to check that the eigenvalues of the partial transpose of $\rho_2$ are $\{1, 0, 0, 0\}.$ As all of these eigenvalues 
are nonnegative this guarantees, by the Peres-Horodecki criterion, \cite{Peres1,Horo}, that $\rho_2$ is indeed a separable state. 

{\bf Example 2:} We consider a 2-qubit quantum state defined by the following density matrix, $\rho_2,$ as given below:
\bes
\ \rho_2
\ =
\bbm
\frac{4}{49}  &  -\frac{6}{35} & \frac{4\sqrt{314}}{735} &  -\frac{4i}{21}  \\
-\frac{6}{35} &  \frac{9}{25} & -\frac{2\sqrt{314}}{175} & \frac{2i}{5}  \\
\frac{4\sqrt{314}}{735}  &  -\frac{2\sqrt{314}}{175} & \frac{1256}{11025} & -\frac{4i\sqrt{314}}{315}  \\
\frac{4i}{21}  &  -\frac{2i}{5} & \frac{4i\sqrt{314}}{315}  & \frac{4}{9} 
\ebm
\ees 
It is easy to see that $tr((\rho_2)^2) = 1$ therefore $\rho_2$ represents a ``pure`` state. 

We now express $\rho_2$ in terms of Generalized Pauli basis and construct the matrix of 
coefficients, $\mathbb{B},$ as given below: 
\bes
\ \mathbb{B}
\ =
\bbm
\frac{1}{4}  &  -\frac{\sqrt{3}}{35} & -\frac{2\sqrt{314}}{315} & -\frac{6713}{4\ts 11025}  \\
\frac{2\sqrt{314}}{735}  &  -\frac{\sqrt{314}}{175} & -\frac{2}{21} & \frac{2\sqrt{314}}{735}  \\
\frac{1}{5}  &  -\frac{2}{21} & -\frac{\sqrt{314}}{175} & -\frac{1}{5}  \\
-\frac{1287}{4\ts 11025}  &  -\frac{3}{35} & -\frac{2\sqrt{314}}{315} & \frac{575}{4\ts 11025} 
\ebm
\ees 
and apply Theorem 4.1 to it. Note that $\alpha_{11} = \frac{1}{4}.$
It is easy to check that the rank of $\mathbb{B}$ is not equal to one therefore the 2-qubit quantum state defined by the 
density matrix, $\rho_2,$ is a {\bf pure entangled} state. 

It is easy to check that the eigenvalues of the partial transpose 
of $\rho_2$ are $\{-0.2780, 0.0844, 0.2780, 0.9156\}.$ As one of these eigenvalues 
is negative this guarantees, by the Peres-Horodecki criterion, \cite{Peres1,Horo}, that $\rho_2$ is indeed an entangled state.

{\bf Example 3:} We consider a 2-qubit quantum state defined by the following density matrix, $\rho_2,$ as given below:
\bes
\ \rho_2
\ =
\bbm
\frac{11}{72}  &  -\frac{i}{18} + \frac{i}{8\sqrt{3}} & \frac{1}{24} & \frac{i}{8\sqrt{3}}  \\
\frac{i}{18} - \frac{i}{8\sqrt{3}}  &  \frac{25}{72} & -\frac{i}{8\sqrt{3}} & \frac{7}{24}  \\
\frac{1}{24}  &  \frac{i}{8\sqrt{3}} & \frac{11}{72} & -\frac{i}{18} + \frac{i}{8\sqrt{3}}  \\
-\frac{i}{8\sqrt{3}}  &  \frac{7}{24} & \frac{i}{18} - \frac{i}{8\sqrt{3}} & \frac{25}{72} 
\ebm
\ees 
It is easy to see that $tr((\rho_2)^2) < 1$ therefore $\rho_2$ represents a ``mixed`` state. 

We now express $\rho_2$ in terms of Generalized Pauli basis and construct the matrix of 
coefficients, $\mathbb{B},$ as given below: 
\bes
\ \mathbb{B}
\ =
\bbm
\frac{1}{4}  &  0 & (\frac{1}{18}-\frac{1}{8\sqrt{3}}) & -\frac{7}{72}  \\
\frac{{1}}{6}  &  0 & -\frac{1}{8\sqrt{3}} & -\frac{1}{8}  \\
0  &  0 & 0 & 0  \\
0  &  0 & 0 & 0 
\ebm
\ees 
As per Algorithm 2 we split the elements of the first row of $\mathbb B$ and rewrite the 
matrix, $\mathbb B,$ as $\mathbb B = \mathbb B^1 + \mathbb B^2$ where
\bes
\ \mathbb{B}^1
\ =
\bbm
\frac{1}{6}  &  0  & -\frac{1}{8\sqrt{3}} & -\frac{1}{8}  \\
\frac{1}{6}  &  0  & -\frac{1}{8\sqrt{3}} & -\frac{1}{8}  \\
      0        &        0       &       0       &       0        \\
      0        &        0       &       0       &       0 
\ebm
\ees 
and 
\bes
\ \mathbb{B}^2
\ =
\bbm
\frac{1}{12}  &  0  & \frac{1}{18} & \frac{1}{36}  \\
      0        &        0       &       0       &       0        \\ 
      0        &        0       &       0       &       0        \\
      0        &        0       &       0       &       0       
\ebm
\ees   
such that now the rank of the matrix $\mathbb{B}^1$ is indeed equal to one (note that the rank of the matrix $\mathbb{B}^2$ 
is obviously equal to one). Therefore, as per Theorem 4.2 $\rho_2$ is a {\bf mixed separable} state. We can now  
easily check in $\mathbb{B}^1$ that (i) $row_2 = 1\ts row_1,$ (ii) $row_3 = 0\ts row_1,$ 
and (iii) $row_4 = 0\ts row_1.$ Therefore, following Algorithm 2 and Theorem 4.2 we have  
$$\rho_2 = \frac{2}{3}(\frac{I}{2}+\frac{{1}}{2}\sigma_1)\ots (\frac{I}{2}-\frac{\sqrt{3}}{8}\sigma_2-\frac{3}{8}\sigma_3) 
+ \frac{1}{3}(\frac{I}{2})\ots (\frac{I}{2}+\frac{1}{3}\sigma_2+\frac{1}{6}\sigma_3)$$ 

It is easy to check that the eigenvalues of the partial transpose of $\rho_2$ are all $\ge 0.$ As all of these eigenvalues 
are nonnegative this guarantees, by the Peres-Horodecki criterion, \cite{Peres1,Horo}, that $\rho_2$ is indeed a separable 
state.  

{\bf Example 4:} We consider a 2-qubit quantum state defined by the following density matrix, $\rho_2,$ as given below:

\bes
\ \rho_2
\ = 
\bbm
0  &  0 & 0 & 0  \\
0  &  \frac{1}{2} & -\frac{1}{2} & 0  \\
0  &  -\frac{1}{2} & \frac{1}{2} & 0  \\
0  &  0 & 0 & 0 
\ebm
\ees 
It is easy to see that $tr((\rho_2)^2) = 1$ therefore $\rho_2$ represents a ``pure'' state. 

We now express $\rho_2$ in terms of Generalized Pauli basis and construct the matrix of 
coefficients, $\mathbb{B},$ as given below: 
\bes
\ \mathbb{B}
\ =
\bbm
\frac{1}{4}  &  0  &  0   &  0 \\
0  &  -\frac{1}{4}  &  0  &  0 \\
0  &  0  &  -\frac{1}{4}  &  0 \\
0  &  0  &  0  &  -\frac{1}{4} 
\ebm
\ees
and apply Theorem 4.1 and Algorithm 1 to it. Note that $\alpha_{11} = \frac{1}{4}.$
It is easy to check that the rank of $\mathbb{B}$ is greater than one. Therefore, the 2-qubit quantum state defined by the 
density matrix, $\rho_2,$ is a {\bf pure entangled} state. 

It is easy to check that the eigenvalues of the partial transpose 
of $\rho_2$ are $\{-\frac{1}{2}, \frac{1}{2}, \frac{1}{2}, \frac{1}{2}\}.$ As one of these eigenvalues 
is negative this guarantees, by the Peres-Horodecki criterion, \cite{Peres1,Horo}, that $\rho_2$ is indeed an entangled state.

{\bf Example 5:} We consider a 2-qubit quantum state defined by the following density matrix, $\rho_2,$ as given below:

\bes
\ \rho_2
\ = 
\bbm
\frac{1}{6}  &  0 & 0 & \frac{1}{6}  \\
0  &  \frac{1}{3} & -\frac{1}{3} & 0  \\
0  &  -\frac{1}{3} & \frac{1}{3} & 0  \\
\frac{1}{6}  &  0 & 0 & \frac{1}{6}
\ebm
\ees 
It is easy to see that $tr((\rho_2)^2) < 1$ therefore $\rho_2$ represents a ``mixed'' state. 

We now express $\rho_2$ in terms of Generalized Pauli basis and construct the matrix of 
coefficients, $\mathbb{B},$ as given below: 
\bes
\ \mathbb{B}
\ =
\bbm
\frac{1}{4}  &  0  &  0   &  0 \\
0  &  -\frac{1}{3}  &  0  &  0 \\
0  &  0  &  -1  &  0 \\
0  &  0  &  0  &  -\frac{1}{3} 
\ebm
\ees
and apply Theorem 4.2 and Algorithm 2 to it. Note that $\alpha_{11} = \frac{1}{4}.$
It is easy to check that $\mathbb{B}$ does not satisfy any of the cases described in Algorithm 2. Therefore, the 2-qubit 
quantum state defined by the density matrix, $\rho_2,$ is a {\bf mixed entangled} state. 

It is easy to check that the eigenvalues of the partial transpose 
of $\rho_2$ are $\{\frac{1}{2}, \frac{1}{2}, -\frac{1}{6}, \frac{1}{6}\}.$ As one of these eigenvalues 
is negative this guarantees, by the Peres-Horodecki criterion, \cite{Peres1,Horo}, that $\rho_2$ is indeed an entangled state. 

{\bf Example 6:} Consider the following 2-qubit {\bf mixed entangled state} defined by the following density matrix, $\rho_2^a.$ 

\bes
\ \rho_2^a
\ = 
\bbm
\frac{1}{3}  &  0 & 0 & \frac{r}{2} \\
0  &  \frac{1}{3} & 0 & 0  \\
0  &  0 & 0 & 0  \\
\frac{r}{2}  &  0 & 0 & \frac{1}{3}
\ebm
\ees
where $0\le r\le \frac{2}{3}.$ This state is a ``maximally'' entangled state \cite{Wei}.

We now express $\rho_2^a$ in terms of Generalized Pauli basis and construct the matrix of 
coefficients, $\mathbb{B}^a,$ as given below: 
\bes
\ \mathbb{B}^a
\ =
\bbm
\frac{1}{4}  &  0  &  0   &  -\frac{1}{12} \\
0  &  \frac{r}{4}  &  0  &  0 \\
0  &  0  &  -\frac{r}{4}  &  0 \\
\frac{1}{12}  &  0  &  0  &  \frac{1}{12} 
\ebm
\ees
and apply Theorem 4.2 and Algorithm 2 to it. Note that $\alpha_{11} = \frac{1}{4}.$
It is easy to check that as long as $r\neq 0$ $\mathbb{B}^a$ does not satisfy any of the cases described in Algorithm 2. Therefore, the 2-qubit 
quantum state defined by the density matrix, $\rho_2^a,$ is a {\bf mixed entangled} state.

It is easy to check that the eigenvalues of the partial transpose 
of $\rho_2^a$ are $\{\frac{1}{6}+(\frac{1}{6}\sqrt{(1+9r^2)}), \frac{1}{6}-(\frac{1}{6}\sqrt{(1+9r^2)}), \frac{1}{3}, \frac{1}{3}\}.$  
The second of these eigenvalues will be negative for $0 < r\le \frac{2}{3}.$ So for such $r$ by the Peres-Horodecki criterion, 
\cite{Peres1,Horo}, $\rho_2^a$ is indeed an entangled state. 

Now, what happens to the state when $r = 0$? We now see below that at $r = 0$ the state becomes a {\bf mixed separable} state.

With $r = 0$ the above density matrix, $\rho_2^a,$ at $r = 0$ becomes 
\bes
\ \rho_2^a
\ = 
\bbm
\frac{1}{3}  &  0 & 0 & 0 \\
0  &  \frac{1}{3} & 0 & 0  \\
0  &  0 & 0 & 0  \\
0  &  0 & 0 & \frac{1}{3}
\ebm
\ees
We now express this $\rho_2^a$ at $r = 0$ in terms of Generalized Pauli basis and construct the matrix of 
coefficients, $\mathbb{B}^a,$ as given below: 
\bes
\ \mathbb{B}^a
\ =
\bbm
\frac{1}{4}  &  0  &  0   &  -\frac{1}{12} \\
0  &  0  &  0  &  0 \\
0  &  0  &  0  &  0 \\
\frac{1}{12}  &  0  &  0  &  \frac{1}{12} 
\ebm
\ees
and apply Theorem 4.2 and Algorithm 2 to it. Note that $\alpha_{11} = \frac{1}{4}.$
It is easy to check that we can split $\mathbb{B}^a$ as $\mathbb{B}^a = \mathbb{B}^{a_1} + \mathbb{B}^{a_2}$ 
where 
\bes
\ \mathbb{B}^{a_1}
\ = 1/3
\bbm
\frac{1}{4}  &  0  &  0   &  \frac{1}{4} \\
0  &  0  &  0  &  0 \\
0  &  0  &  0  &  0 \\
\frac{1}{4}  &  0  &  0  &  \frac{1}{4} 
\ebm
\ees
and 
\bes
\ \mathbb{B}^{a_2}
\ = 2/3
\bbm
\frac{1}{4}  &  0  &  0   &  -\frac{1}{4} \\
0  &  0  &  0  &  0 \\
0  &  0  &  0  &  0 \\
0  &  0  &  0  &  0 
\ebm
\ees
Using above splitting of $\mathbb{B}^a$ we can easily write $\rho_2^a$ at $r = 0$ as follows:

$$\rho_2^a = \frac{1}{3}[(\frac{I}{2} + \frac{\sigma_3}{2})\ots (\frac{I}{2} + \frac{\sigma_3}{2})] + \frac{2}{3}[(\frac{I}{2})\ots (\frac{I}{2}-\frac{\sigma_3}{2})]$$
Therefore, when $r = 0$ the 2-qubit 
quantum state defined by the density matrix, $\rho_2^a,$ becomes a {\bf mixed separable} state.

It is easy to check that the eigenvalues of the partial transpose 
of $\rho_2^a$ at $r = 0$ are $\{\frac{1}{3},\frac{1}{3},\frac{1}{3},\frac{1}{3}\}.$  
As all of these eigenvalues are now positive (nonnegative) therefore by the Peres-Horodecki criterion, 
\cite{Peres1,Horo}, $\rho_2^a$ for $r = 0$ case is indeed a separable state. 

For the other range, namely, $\frac{2}{3}\le r\le 1,$ we have 
\bes
\ \rho_2^b
\ = 
\bbm
\frac{r}{2}  &  0 & 0 & \frac{r}{2} \\
0  &  \frac{1}{3} & 0 & 0  \\
0  &  0 & 0 & 0  \\
\frac{r}{2}  &  0 & 0 & \frac{1}{3}
\ebm
\ees
We now express $\rho_2^b$ in terms of Generalized Pauli basis and construct the matrix of 
coefficients, $\mathbb{B}^b,$ as given below: 
\bes
\ \mathbb{B}^b
\ =
\bbm
(\frac{r}{8}+\frac{1}{6})  &  0  &  0   &  (\frac{r}{8}-\frac{1}{6}) \\
0  &  \frac{r}{4}  &  0  &  0 \\
0  &  0  &  -\frac{r}{4}  &  0 \\
\frac{r}{8}  &  0  &  0  &  \frac{r}{8} 
\ebm
\ees
and apply Theorem 4.2 and Algorithm 2 to it. 
It is easy to check that for the entirs range of $r$ $\mathbb{B}^b$ does not satisfy any of the cases described in Algorithm 2. Therefore, the 2-qubit 
quantum state defined by the density matrix, $\rho_2^b,$ is a {\bf mixed entangled} state.

{\bf Example 7:} In this example and the next example we will see that the density matrices apparently look 
similar in the form as they both contain the {\bf same} $2\ts 2$ matrix
\bes
\bbm
\frac{1}{2} & \frac{1}{2} \\
\frac{1}{2} & \frac{1}{2}
\ebm
\ees
surrounded by zeros. But we will see below that the density matrix in this example corresponds to a 2-qubit ``entangled'' state 
while the density matrix in the next example corresponds to a 3-qubit ``separable'' state. The 2-qubit state we consider here 
is defined by the following density matrix, $\rho_2:$ 
\bes
\ \rho_2
\ = 
\bbm
0  &  0 & 0 & 0  \\
0  &  \frac{1}{2} & \frac{1}{2} & 0  \\
0  &  \frac{1}{2} & \frac{1}{2} & 0  \\
0  &  0 & 0 & 0 
\ebm
\ees 
It is easy to see that $tr((\rho_2)^2) = 1$ therefore $\rho_2$ represents a ``pure'' state. 

We now express $\rho_2$ in terms of Generalized Pauli basis and construct the matrix of 
coefficients, $\mathbb{B},$ as given below: 
\bes
\ \mathbb{B}
\ =
\bbm
\frac{1}{4}  &  0  &  0   &  0 \\
0  &  \frac{1}{4}  &  0  &  0 \\
0  &  0  &  \frac{1}{4}  &  0 \\
0  &  0  &  0  &  -\frac{1}{4} 
\ebm
\ees
and apply Theorem 4.1 and Algorithm 1 to it. Note that $\alpha_{11} = \frac{1}{4}.$
It is easy to check that the rank of $\mathbb{B}$ is greater than one. Therefore, the 2-qubit quantum state defined by the 
density matrix, $\rho_2,$ is a {\bf pure entangled} state. 

{\bf Example 8:} We now consider a similar looking 3-qubit quantum state defined by the following density matrix, $\rho_3,$ 
which as we will see below is a {\bf pure separable} state:
\bes
\ \rho_3
\ = 
\bbm
0  &  0 & 0 & 0 & 0 & 0 & 0 & 0 \\
0  &  0 & 0 & 0 & 0 & 0 & 0 & 0 \\
0  &  0 & 0 & 0 & 0 & 0 & 0 & 0 \\
0  &  0 & 0 & 0 & 0 & 0 & 0 & 0 \\
0  &  0 & 0 & 0 & \frac{1}{2} & \frac{1}{2} & 0 & 0 \\
0  &  0 & 0 & 0 & \frac{1}{2} & \frac{1}{2} & 0 & 0 \\
0  &  0 & 0 & 0 & 0 & 0 & 0 & 0 \\
0  &  0 & 0 & 0 & 0 & 0 & 0 & 0 \\ 
\ebm
\ees
It is easy to see that $tr((\rho_2)^2) = 1$ therefore $\rho_3$ represents a ``pure'' state. Therefore we proceed to apply 
Theorem 5.1 and Algorithm 3. It is easy to check that 
\bes
\ A_0
\ = 
\bbm
\frac{1}{8}  & \frac{1}{8}  & 0 & 0  \\
0  &  0 & 0 & 0  \\
0  &  0 & 0 & 0  \\
\frac{1}{8}  &  \frac{1}{8} & 0 & 0 
\ebm
\ees 
Therefore the rank of $A_0$ is equal to one. 

It is easy to check further that $A_1 = 0\ts A_0,$ $A_2 = 0\ts A_0,$ and $A_3 = (-1)\ts A_0.$ Therefore, as per Algorithm 3
the state defined by the above given density matrix $\rho_3$ is a {\bf pure separable} state. We can easily find the following 
expression for this density matrix, $\rho_3,$ consistent with Definition 3.4, by proceeding as per Algorithm 3, Step 2 (a):

$\rho_3 = (\rho_1^1)\ots (\rho_1^2)\ots (\rho_1^3) = (\frac{I}{2}-\frac{\sigma_3}{2})\ots (\frac{I}{2}+\frac{\sigma_3}{2})\ots (\frac{I}{2}+\frac{\sigma_1}{2}).$

{\bf Example 9:} We consider a 3-qubit quantum state defined by the following density matrix, $\rho_3,$ as given below:

\bes
\ \rho_3
\ = 
\bbm
\frac{1}{4}  &  -\frac{1}{4} & 0 & 0 & 0 & 0 & \frac{1}{4} & -\frac{1}{4} \\
-\frac{1}{4}  &  \frac{1}{4} & 0 & 0 & 0 & 0 & -\frac{1}{4} & \frac{1}{4} \\
0  &  0 & 0 & 0 & 0 & 0 & 0 & 0 \\
0  &  0 & 0 & 0 & 0 & 0 & 0 & 0 \\
0  &  0 & 0 & 0 & 0 & 0 & 0 & 0 \\
0  &  0 & 0 & 0 & 0 & 0 & 0 & 0 \\
\frac{1}{4}  &  -\frac{1}{4} & 0 & 0 & 0 & 0 & \frac{1}{4} & -\frac{1}{4} \\
-\frac{1}{4}  &  \frac{1}{4} & 0 & 0 & 0 & 0 & -\frac{1}{4} & \frac{1}{4} 
\ebm
\ees
It is easy to see that $tr((\rho_2)^2) = 1$ therefore $\rho_3$ represents a ``pure'' state. Therefore we proceed to apply 
Theorem 5.1 and Algorithm 3. It is easy to check that 
\bes
\ A_0
\ = 
\bbm
\frac{1}{8}  & -\frac{1}{8}  & 0 & 0  \\
0  &  0 & 0 & 0  \\
0  &  0 & 0 & 0  \\
\frac{1}{8}  &  -\frac{1}{8} & 0 & 0 
\ebm
\ees 
Therefore the rank of $A_0$ is equal to one. But
\bes
\ A_1
\ = 
\bbm
0  & 0  & 0 & 0  \\
\frac{1}{8}  &  -\frac{1}{8} & 0 & 0  \\
0  &  0 & 0 & 0  \\
0  &  0 & 0 & 0 
\ebm
\ees
therefore, $A_1 \neq k_1\ts A_0,$ where $k_1$ a constant, and so $\rho_3$ is a {\bf pure entangled} state as per (c) in Algorithm 3.

{\bf Remarks:} We can proceed as above and can work out some examples for 3-qubit and 4-qubit cases. We have not done it here 
since it involves huge calculations. We just make few remarks in this regard and stop. For 3-qubit and 4-qubit cases 
the density matrices will be of sizes $8\ts 8$ and $16\ts 16$ respectively. 
The matrices that forms the Generalized Pauli basis for 3-qubit and 4-qubit cases will also be of sizes $8\ts 8$ and $16\ts 16$ 
respectively. The cardinality of the basis vectors for 3-qubit and 4-qubit cases (in terms of the matrices of sizes $8\ts 8$ 
and $16\ts 16$ respectively) will be 64 and 256 in number respectively. Thus, the size of the computation involved for solving problems in the 
higher dimensional cases increases very rapidly although the procedures for solving these problems remain quite similar. 

\section{A New Quantum Protocol for Superluminal Classical Communication:} 

Entanglement is a resource of great utility in quantum computation and 
quantum information \cite{nie,Colin}. Schrödinger described it as follows: ``I would not call [entanglement] one but 
rather the characteristic trait of quantum mechanics, the one that enforces its entire departure from classical lines of 
thought.'' \cite{Sch}. Being usually fragile to environment, entanglement is robust against conceptual and mathematical 
tools that struggle with the task of deciphering its rich structure \cite{Horo}. Entanglement is a fundamental resource 
of nature which is comparable in importance to energy, information, entropy, or any other equally fundamental resource. 
Quantum entanglement describes the strong correlation that exists among different parts of composite quantum systems even 
after the parts of these composite quantum systems may get spacelike separated from one another. The so called teleportation 
of any arbitrary quantum state  from one party to other party \cite{CB} has been verified experimentally. This quantum 
teleportation is not only a theoretical idea but has been verified in the teleportation 
experiments between two locations separated over the distance of the order of hundred kilometers \cite{Bou,Yin,Ma,Lan,Sun,Val}. 
Recently a team of Chinese scientists have extended this range to about 1400 kilometers \cite{Chinese}. 

In this paper we propose a new quantum protocol to instantaneously and simultaneously transmit any desired classical information 
in terms of any desired sequence of classical bits from Alice stationed at a data transmission center to all the 
other recipients, located at different faraway locations from Alice at the data transmission center. These 
participants are situated at the locations which are faraway not only from Alice but also could be faraway from one  
another too. Thus, the proposed quantum protocol makes use of the entangled shared states and the nonlocal influences to transmit 
the desired known classical information in terms of a sequence of classical bits, simultaneously and instantaneously, 
from Alice sitting at the data transmission center to different distant locations, faraway from data transmission center 
and could be faraway from one another too.  

We will also see that the new quantum protocol proposed here for the classical data transmission can be used as an aid in 
the well-known quantum teleportation protocol \cite{CB}. 
The quantum teleportation protocol \cite{CB} requires to transmit two classical bits from Alice to Bob 
to complete the teleportation protocol successfully. This was so far impossible to do faster than the velocity of light. 
But the new quantum protocol proposed here the classical data transmission can achieve the classical data transmission in 
terms of a data sequence 
of classical bits instantaneously. This makes possible the instantaneous transmission of two classical bits from Alice to Bob 
for Bob to successfully complete the teleportation protocol \cite{CB}.
In short, the proposed quantum protocol will work as a superluminal communicator to pass two classical bits 
that appear during the Bell basis measurement by Alice to the now faraway Bob to perform required operation, and as a 
consequence makes possible the instantaneous teleportation of a quantum state from Alice to Bob. 

\section{Some Preliminaries for a New Quantum Protocol for Superluminal Classical Communication:}  

We now give some preliminaries required for the proposed quantum protocol for superluminal classical communication. 

Among the 2-qubit states the well-known four Bell basis states, $|\beta_{ij}\ran, i,j \in \{0,1\},$ are the maximally 
entangled states. They are defined as $$|\beta_{00}\ran = \frac{1}{\sqrt{2}}(|00\ran + |11\ran)$$ $$|\beta_{01}\ran = \frac{1}{\sqrt{2}}(|01\ran + |10\ran)$$ 
$$|\beta_{10}\ran = \frac{1}{\sqrt{2}}(|00\ran - |11\ran)$$ $$|\beta_{11}\ran = \frac{1}{\sqrt{2}}(|01\ran - |10\ran)$$

Both these sets of states, the (2-qubit) computational basis states and the (2-qubit) Bell basis states, form orthonormal 
basis so that any 2-qubit pure quantum state can be expressed in terms of them. 

Also, the four Bell basis states, $|\beta_{ij}\ran, i,j \in \{0,1\},$ and the four computational basis states, 
$|ij\ran, i,j \in \{0,1\}$ are 
inter convertible into each other through the following invertible matrix, $A$ , which is also unitary, i.e. $A^{-1} = A^{+}$ 
where the matrix $A^{+}$ denotes the conjugate transpose of the matrix $A.$

The matrix $A$ is as given below:

\bes
\ A
\ =
\frac{1}{\sqrt{2}}
\bbm
 1 &  0 & 0 & 1 \\
0  &  1 & 1 & 0 \\
1 &  0 & 0 & -1 \\
0  &  1 & -1 & 0 
\ebm
\ees

and the Bell basis and the computational basis are inter convertible and are related by the following matrix equations:

\bes
\bbm
|\beta_{00}\ran \\
|\beta_{01}\ran \\
|\beta_{10}\ran \\
|\beta_{11}\ran 
\ebm
\ =
A
\bbm
|00\ran \\
|01\ran \\
|10\ran \\
|11\ran 
\ebm
\ees

and

\bes
\bbm
|00\ran \\
|01\ran \\
|10\ran \\
|11\ran 
\ebm
\ =
A^{+}
\bbm
|\beta_{00}\ran \\
|\beta_{01}\ran \\
|\beta_{10}\ran \\
|\beta_{11}\ran 
\ebm
\ees

One can easily extend the idea and define the Bell basis states for the 3-qubit states. For the 3-qubit case we will have 
eight Bell basis states $|\beta_{ijk}\ran, i,j,k \in \{0,1\}$ and they are related to eight computational basis 
states, $|ijk\ran, i,j,k \in \{0,1\}$ through the following matrix equation:

\bes
\bbm
|\beta_{000}\ran \\
|\beta_{001}\ran \\
|\beta_{010}\ran \\
|\beta_{011}\ran \\
|\beta_{100}\ran \\
|\beta_{101}\ran \\
|\beta_{110}\ran \\
|\beta_{111}\ran 
\ebm
\ =
B
\bbm
|000\ran \\
|001\ran \\
|010\ran \\
|011\ran \\
|100\ran \\
|101\ran \\
|110\ran \\
|111\ran 
\ebm
\ees

where the matrix $B$ is as given below:

\bes
\ B
\ =
\frac{1}{\sqrt{2}}
\bbm
1 &  0 & 0 &  0 & 0 &  0 & 0 & 1 \\
0 &  1 & 0 &  0 & 0 &  0 & 1 & 0 \\
0 &  0 & 1 &  0 & 0 &  1 & 0 & 0 \\
0 &  0 & 0 &  1 & 1 &  0 & 0 & 0 \\
1 &  0 & 0 &  0 & 0 &  0 & 0 & -1 \\
0 &  1 & 0 &  0 & 0 &  0 & -1 & 0 \\
0 &  0 & 1 &  0 & 0 &  -1 & 0 & 0 \\
0 &  0 & 0 &  1 & -1 &  0 & 0 & 0
\ebm
\ees

It is easy to check that like the matrix $A$ defined above the $8\ts 8$ matrix $B$ is also invertible and unitary, 
i.e.$B^{-1}$ exists and $B^{-1} = B^{+}$ 
where the matrix $B^{+}$ denotes the conjugate transpose of the matrix $B.$ 

A further generalization is also possible. We can proceed 
on the similar lines and define orthonormal Bell basis containing $2^n$ Bell basis states for $n$-qubit case, 
namely,$|\beta_{i_1i_2\cdots i_n}\ran, i_1,i_2,\cdots, i_n \in \{0,1\}$  and relate them with the corresponding 
$n$-qubit computational basis states, namely, $|i_1i_2\cdots i_n\ran$ such that $i_1,i_2,\cdots, i_n \in \{0,1\},$ 
through a $2^n\ts 2^n$ matrix which will be also invertible and unitary.

A formula for writing an $n$-qubit Bell basis state in terms of $n$-qubit computational basis states is as given below: 

$$|\beta_{i_1i_2\cdots i_n}\ran = \frac{1}{\sqrt{2}}(|0i_2i_3\cdots i_n\ran + (-1)^{i_1}|1(1-i_2)(1-i_3)\cdots (1-i_n)\ran)$$

Using this relation one can easily construct the transformation matrix of size $2^n\ts 2^n$ to express $n$-qubit Bell basis 
states in terms of corresponding $n$-qubit computational basis states. One can further see that as in the previous cases this 
matrix will also be invertible and unitary. 

One can easily build the quantum circuit to obtain the $n$-qubit Bell basis state. By operating this quantum circuit on the 
computational basis state, $|x_1x_2\cdots x_n\ran,$ one will obtain the desired Bell basis state, 
$|\beta_{x_1x_2\cdots x_n}\ran.$ 

One can build this quantum circuit by carrying out the following steps:  

(i) INPUT: Start with the computational basis state, $|x_1x_2\cdots x_n\ran.$ 

(ii) Operate the Hadamard gate, $H,$ on the first qubit $|x_1\ran.$ 

(iii) Operate the CNOT gates, $C_{12},C_{13},C_{14},\cdots, C_{1n}$ with the first qubit, namely, $H|x_1\ran$ as the control 
qubit and respectively $|x_2\ran,|x_3\ran,|x_4\ran,\cdots, |x_n\ran$ as the target qubits. 

(iv) OUTPUT: One arrives at the desired Bell basis state, $|\beta_{x_1x_2\cdots x_n}\ran.$

It is interesting to note that the computational basis state, $|x_1x_2\cdots x_n\ran,$ is a separable state while the Bell 
basis state, $|\beta_{x_1x_2\cdots x_n}\ran,$ is an entangled state, in fact maximally entangled state.

\section{A Quantum Protocol for Transmitting ``Same'' Information to multiple Number of Recipients:} 

In this section we discuss the new quantum protocol for transmitting the information in the form of a sequence of classical 
bits (zeros and ones) arranged in a particular order from Alice at data transmission center to all the recipients located 
at faraway locations.  
There are several recipients situated faraway from data transmission center and also from each other and suppose 
we want to transmit the ``identical'' information to all, bit by bit, in terms of classical bits, instantaneously and 
simultaneously. 

Is it possible to fulfill this task successfully? 

The answer to this question is YES. In this section we will describe how the new quantum protocol developed here 
to transmit the classical information can pass this classical information, not only correctly but also instantaneously and simultaneously, 
in terms of classical bits from Alice to all the recipients located at different distant locations. 

For the sake of simplicity we first consider here the simple case where Alice stationed at data transmission center will 
manage to transmit, bit by bit, the ``same'' data sequence $s_1,s_2, \cdots, s_m$ where $s_i\in \{0,1\},$ to all the recipients 
say $P_1,P_2,P_3,\cdots,P_n$ which are faraway not only from Alice but from each other too. 

The proposed quantum protocol proceeds as follows: 

Initially, Alice and all the participants are together at the data transmission center and when they are together at the 
data transmission center they prepare several copies of the following two types of ``Generalized'' Bell basis states 
$|\Delta_1\ran$ and $|\Delta_2\ran$:

$$|\Delta_1\ran = |\beta_{0(j_1j_2\cdots j_r)(j_1j_2\cdots j_r)\cdots (j_1j_2\cdots j_r)}\ran,$$ shared 
by Alice and all the other participants of the protocol say $P_1,P_2,P_3,\cdots,P_n.$ The first qubit of each of 
these ``Generalized'' Bell basis states is with Alice, the next group of r qubits represented above 
as $(j_1j_2\cdots j_r)$ is the first group of r qubits and it is in the possession of participant $P_1$ 
(as a ket $|j_1j_2\cdots j_r\ran.$) 
The next group of r qubits represented above as $(j_1j_2\cdots j_r)$ and it is in the possession of participant $P_2$ 
(as a ket $|j_1j_2\cdots j_r\ran,$) and so on. 
The last group of last r qubits represented above as $(j_1j_2\cdots j_r)$ and it is in the possession of participant $P_n$ 
(as a ket $|j_1j_2\cdots j_r\ran.$) Similarly, 
$$|\Delta_2\ran = |\beta_{0(k_1k_2\cdots k_r)(k_1k_2\cdots k_r)\cdots (k_1k_2\cdots k_r)}\ran,$$ shared 
by Alice and all the other participants of the protocol say $P_1,P_2,P_3,\cdots,P_n.$ The first qubit of each of 
these ``Generalized'' Bell basis states is with Alice, the group of next r qubits, as above,  
represented as $(k_1k_2\cdots k_r)$ is the first group of r qubits and it is in the possession of participant $P_1$ 
(as a ket $|k_1k_2\cdots k_r\ran.$) 
The next group of r qubits represented above as $(k_1k_2\cdots k_r)$ and it is in the possession of participant $P_2$ 
(as a ket $|k_1k_2\cdots k_r\ran,$) and so on. 
The last group of last r qubits represented above as $(k_1k_2\cdots k_r)$ and it is in the possession of participant $P_n$ 
(as a ket $|k_1k_2\cdots k_r\ran.$) 

Note that when Alice and all the participants were together at the data transmission center they prepared several copies 
of the above mentioned ``Generalized'' Bell basis states, $|\Delta_1\ran$ and $|\Delta_2\ran$. 

All the participants $P_1,P_2,P_3,\cdots,P_n$  
who will be the recipients of information transmitted from the data transmission center by Alice then go away to 
their respective faraway locations along with their qubits while Alice always stays at the data transmission center.

Alice is then provided with the classical information in terms of a sequence of classical bits 
$s_1,s_2, \cdots, s_m$ where $s_i\in \{0,1\}.$ This sequence of classical bits is to be transmitted bit by bit to all the 
faraway recipients in the same order as given, i.e. firstly the value of the classical bit $s_1,$ which is given to Alice, 
to transmit to all the participants, $P_1,P_2,P_3,\cdots,P_n,$ 
who have now gone faraway from Alice as well as from each other, instantaneously and simultaneously. Secondly, Alice is given 
$s_2$ to transmit to all participants, instantaneously and simultaneously, and so on, till finally Alice has to transmit 
$s_m$ to all the participants, instantaneously and simultaneously.

Note that the main task for Alice is to manage to transmit this information in terms of a known sequence of classical bits 
$s_1,s_2, \cdots, s_m$ where $s_i\in \{0,1\},$ such that it is transmitted to all the faraway recipients and each bit 
transmitted by Alice must reach all the recipients instantaneously and simultaneously.

We now explain what information we have associated with the following sequences of classical 
bits (zeros and ones), namely, (i) $j_1j_2\cdots j_r$ and (ii) $k_1k_2\cdots k_r$. We have chosen these sequences such that 
the sequence ``$j_1j_2\cdots j_r$'' stands for the information, say ``BLUE SKY'', and the sequence 
``$(1-k_1)(1-k_2)\cdots (1-k_r)$'' stands for the information, say ``BRIGHT SUN''. Similarly, the information associated by us 
with the following sequences ``$(1-j_1)(1-j_2)\cdots (1-j_r)$'' and ``$k_1k_2\cdots k_r$'' is ``COMPLEMENT of BLUE SKY'' 
and ``COMPLEMENT of BRIGHT SUN'' respectively. Note that 

``COMPLEMENT of the COMPLEMENT = ORIGINAL''.

Note that Alice will receive the information in terms of an ordered sequence of zeros and ones after all other participants 
will go away to their respective faraway locations and Alice wants to convey the same information just now received to her  
i.e. the sequence of classical bits $s_1,s_2, \cdots, s_m$ say, where $s_i\in \{0,1\}.$ This information is to be conveyed by 
Alice to all the faraway recipients, instantaneously and simultaneously. 

Suppose $s_1 = 0.$ In this case Alice will measure partially the state $|\Delta_1\ran$ among the several copies available to 
her in computational basis. As a result of this measurement the state $|\Delta_1\ran$ will collapse either into state 

$$|\psi_1\ran = |0(j_1j_2\cdots j_r)(j_1j_2\cdots j_r)\cdots (j_1j_2\cdots j_r)\ran$$ or into state 

$$|\psi_2\ran = |1((1-j_1)(1-j_2)\cdots (1-j_r))((1-j_1)(1-j_2)\cdots (1-j_r))\cdots (1-j_1)(1-j_2)\cdots (1-j_r))\ran$$

As a consequence of this partial measurement by Alice in the computational basis she will find the single qubit in her 
possession to be collapsed into either $|0\ran$ or $|1\ran,$ and correspondingly all the remote participants will 
see their state of $r$ qubits in their possession collapsed into either $|j_1j_2\cdots j_r\ran$ 
or $|(1-j_1)(1-j_2)\cdots (1-j_r)\ran$ which conveys the information ``BLUE SKY'' or its COMPLEMENT to all the remote 
participants and it is already told to all the remote participants (when 
they all were together with Alice) that the information ``BLUE SKY'' or its COMPLEMENT corresponds to $s_1 = 0.$ In other 
words, a CONVENTION is made prior to the separation from each other of Alice and all other participants, when they were 
together that the outcome of measurement equal to ``BLUE SKY'' or its COMPLEMENT corresponds to $s_1 = 0.$ 
Note that the COMPLEMENT may be something meaningful arrangement of letters or it could be be some meaningless arrangement 
of letters in the language. 

Note that we have made a convention that whatever outcome we may get, namely, ``BLUE SKY'' or its either meaningful or 
meaningless COMPLEMENT, we have decided a priory that such outcome corresponds to $s_1 = 0.$ Thus Alice 
has managed to convey to all the remote participants that $s_1=0$ by partially measuring in the computational basis a copy of 
the state $|\Delta_1\ran$ among the several copies available to her. 

Suppose $s_1 = 1.$ In this case Alice will measure the first qubit of the state state $|\Delta_2\ran$ in her possession, 
i.e. measure partially the first qubit in her possession of the state $|\Delta_2\ran$ among the several copies available to 
her in the computational basis. As a result of this measurement the state $|\Delta_2\ran$ will collapse either into state 

$$|\phi_1\ran = |0(k_1k_2\cdots k_r)(k_1k_2\cdots k_r)\cdots (k_1k_2\cdots k_r)\ran$$ or into state 

$$|\phi_2\ran = |1((1-k_1)(1-k_2)\cdots (1-k_r))((1-k_1)(1-k_2)\cdots (1-k_r))\cdots ((1-k_1)(1-k_2)\cdots (1-k_r))\ran$$

Now, as a consequence Alice will see her qubit as $|0\ran$ or $|1\ran$ and correspondingly all the remote participants will 
see their state as either $|k_1k_2\cdots k_r\ran$ or $|(1-k_1)(1-k_2)\cdots (1-k_r)\ran$ which conveys the information: a 
``COMPLEMENT of BRIGHT SUN'' or ``BRIGHT SUN'' to all the remote participants and it is already told to all the remote 
participants, when they all were together with Alice, that the information ``BRIGHT SUN'' or its COMPLEMENT corresponds 
to $s_1 = 1.$ In other 
words, a CONVENTION is made prior to the separation from each other of Alice and all other participants, when they were 
together that the outcome of measurement equal to ``COMPLEMENT of BRIGHT SUN'' or ``BRIGHT SUN'' corresponds to $s_1 = 1.$ 
Thus Alice has thus managed to convey to all the remote participants that $s_1=1$ by partially measuring 
in the computational basis a copy of the state $|\Delta_2\ran$ among the several copies available to her. 

By repeating these actions one can then transmit $s_2,s_3, \cdots, s_m$ where $s_i\in \{0,1\},$ and thus the entire sequence 
can be transmitted, bit by bit, to all the recipients, instantaneously and simultaneously. 

\section{A Quantum Protocol for Superluminal Teleportation:}

The quantum teleportation protocol \cite{CB} involves two parties, Alice and Bob. Alice starts with a single qubit 
quantum state, $|\psi\ran = a|0\ran + b|1\ran,$ in her possession which is unknown to her except that 
$|a|^2 + |b|^2 = 1.$ Alice and Bob also start out by sharing between them a Bell state, 
$|\beta_{00}\ran = \frac{1}{\sqrt{2}}[|00\ran + |11\ran]$ such that the 
first qubit of this Bell state is in the possession of Alice and the second qubit of this Bell state is with Bob.  
Suppose that Alice wants to teleport the single qubit state $|\psi\ran$ in her possession to Bob. 

The joint state of three qubits is 
$$|\psi\ran|\beta_{00}\ran = \frac{1}{\sqrt{2}}[a|000\ran + a|011\ran + b|100\ran + b|111\ran].$$ 
Rewriting the above quantum state we have  
$$|\psi\ran|\beta_{00}\ran = \frac{1}{\sqrt{2}}[|00\ran a|0\ran + |01\ran a|1\ran + |10\ran b|0\ran + |11\ran b|1\ran].$$ 
By expressing the 
computational basis states made up of first two qubits $|00\ran, |01\ran, |10\ran, |11\ran$ in the above equation in terms of 
the standard Bell basis states $|\beta_{00}\ran, |\beta_{01}\ran, |\beta_{10}\ran, |\beta_{11}\ran$ where 
$|\beta_{00}\ran = \frac{1}{\sqrt{2}}[|00\ran + |11\ran],$ $|\beta_{01}\ran = \frac{1}{\sqrt{2}}[|01\ran + |10\ran],$ 
$|\beta_{10}\ran = \frac{1}{\sqrt{2}}[|00\ran - |11\ran],$ and $|\beta_{11}\ran = \frac{1}{\sqrt{2}}[|01\ran - |10\ran]$
we get
$$|\psi\ran|\beta_{00}\ran = \frac{1}{2}[|\beta_{00}\ran(a|0\ran + b|1\ran) +  |\beta_{01}\ran(a|1\ran + b|0\ran) + 
|\beta_{10}\ran(a|0\ran - b|1\ran) + |\beta_{11}\ran(a|1\ran - b|0\ran)].$$ It is easy to check further that the above 
equation can also be expressed as follows:
$$|\psi\ran|\beta_{00}\ran = \frac{1}{2}[|\beta_{00}\ran(I|\psi\ran) +  |\beta_{01}\ran(X|\psi\ran) + 
|\beta_{10}\ran(Z|\psi\ran) + |\beta_{11}\ran(X.Z|\psi\ran)],$$ where $I, X, Z$ are standard Pauli operators.

Now, if Alice will perform the partial measurement on the first two qubits in her possession, i.e. if she will perform Bell 
basis measurement on the two qubits in her possession then any one out of the four Bell basis states 
$|\beta_{00}\ran, |\beta_{01}\ran, |\beta_{10}\ran, |\beta_{11}\ran$ will be the outcome with equal probability  
equal to $\frac{1}{4}.$ After Bell basis measurement by Alice yielding one out of the four Bell basis states what  
the posteriory state (third qubit) Bob will have will depend on the outcome of Alice's Bell basis measurement. So, the next 
step of the well known protocol \cite{CB} is to convey the result of the Bell basis measurement done by Alice to Bob 
in terms of two classical bits over a classical channel so that Bob can perform the appropriate recovery operation for 
yielding the exactly identical copy of $|\psi\ran$ as his qubit (third qubit). Now, (without the knowledge of the quantum 
protocol for instantaneously transmitting the classical information, which we have developed in this paper) conveying the 
result of the Bell basis measurement in terms of two classical bits from Alice to Bob over a classical channel cannot be 
done faster than the velocity of light which is the well-known experimentally verified upper limit to the velocity in the 
universe. 

In order to manage the above said classical communication in terms of two classical bits say $s_1s_2$, generated as 
a result of the Bell basis measurement done by Alice, instantaneously, the following ``Generalized'' Bell basis states 
will be prepared when Alice and Bob were together. When together, Alice and Bob will prepare sufficiently many copies of the 
following ``Generalized'' Bell basis states, $|\alpha_1\ran$ and $|\alpha_2\ran$:
$$|\alpha_1\ran = |\beta_{0j_1j_2\cdots j_r}\ran$$ The first qubit of each of the copy of  
these ``Generalized'' Bell basis states, $|\alpha_1\ran,$ is with Alice, and the group of next r qubits represented above 
as $j_1j_2\cdots j_r$ in each of that copy is in the possession of Bob. 
$$|\alpha_2\ran = |\beta_{0k_1k_2\cdots k_r}\ran$$ As above the first qubit of each of the copy of 
these ``Generalized'' Bell basis states, $|\alpha_2\ran,$ is with Alice, and the group of next r qubits represented above 
as $k_1k_2\cdots k_r$ in each of that copy is in the possession of Bob. 

Now, as in done in section III We have chosen these sequences such that 
the sequence ``$j_1j_2\cdots j_r$'' stands for the information, say ``BLUE SKY'', 
and the sequence ``$(1-k_1)(1-k_2)\cdots (1-k_r)$'' stands for the information, say ``BRIGHT SUN''. 

Suppose $s_1 = 0.$ In this case Alice will measure partially the state $|\alpha_1\ran$ among the several copies available to 
her in the computational basis. As a result of this measurement the state $|\alpha_1\ran$ will collapse either into state 

$$|\mu_1\ran = |0j_1j_2\cdots j_r\ran$$ or into state 

$$|\mu_2\ran = |1(1-j_1)(1-j_2)\cdots (1-j_r)\ran$$

As a consequence of this partial measurement by Alice in the computational basis she will find the single qubit in her 
possession to be collapsed into either $|0\ran$ or $|1\ran,$ state and correspondingly Bob will 
see his state of $r$ qubits in his possession collapsed into either $|j_1j_2\cdots j_r\ran$ 
or $|(1-j_1)(1-j_2)\cdots (1-j_r)\ran$ which conveys the information ``BLUE SKY'' or the ``COMPLEMENT of BLUE SKY'' to 
Bob and it is already told to Bob, when 
Alice and Bob were together, that the information ``BLUE SKY'' or the ``COMPLEMENT of BLUE SKY'' corresponds to $s_1 = 0.$ 
Note that the COMPLEMENT may be by chance a meaningful arrangement of letters or it could be be some meaningless arrangement 
of letters in the language.  

Note that we have made a convention that whatever outcome we may get, namely, ``BLUE SKY'' or its either meaningless or 
meaningful COMPLEMENT, we have decided a priory that such outcome corresponds to $s_1 = 0.$ Thus Alice 
has managed to convey to Bob that $s_1=0$ by partially measuring in the computational basis a copy of 
the state $|\alpha_1\ran$ among the sufficiently many copies available to her. 

Suppose $s_1 = 1.$ In this case Alice will measure partially the state $|\alpha_2\ran$ among the several copies available to 
her in the computational basis. As a result of this measurement the state $|\alpha_2\ran$ will collapse either into state 

$$|\nu_1\ran = |0k_1k_2\cdots k_r\ran$$ or into state 

$$|\nu_2\ran = |1(1-k_1)(1-k_2)\cdots (1-k_r)\ran$$
As a consequence of this partial measurement Alice will see her qubit as $|0\ran$ or $|1\ran$ and correspondingly Bob will 
see his state as either $|k_1k_2\cdots k_r\ran$ or $|(1-k_1)(1-k_2)\cdots (1-k_r)\ran$ which conveys the information: a 
``COMPLEMENT of BRIGHT SUN'' or ``BRIGHT SUN'' to Bob and it is already told to Bob,  
when Alice and Bob were together, that the information ``COMPLEMENT of BRIGHT SUN'' or ``BRIGHT SUN'' corresponds 
to $s_1 = 1.$ Thus Alice has managed to convey to Bob that $s_1=1$ by partially measuring 
in the computational basis a copy of the state $|\alpha_2\ran$ among the sufficiently many copies available to her. 

By repeating these actions Alice can also transmit the value of $s_2$ to Bob and thus the classical bits, $s_1,s_2$ that 
appeared due to Bell basis measurement by Alice can be transmitted to Bob instantaneously and using these bits Bob 
will manage to change the (third) qubit in his possession into the desired state, $|\psi\ran.$ The desired state, 
$|\psi\ran,$ which was initially with Alice and got destroyed during her Bell basis measurement, will be made to reincarnate 
by Bob after instantaneously receiving the desired classical bits, $s_1,s_2$ using our new quantum protocol. 
  
{\bf Remarks:}

The ideas developed in the sections IX onwards provide a new technique for instantaneously and simultaneously transmitting 
any piece of known classical information to faraway locations in the universe and thus provide us a means to remain connected 
through instantaneous and simultaneous exchange of information to distant parts of the universe. In short this new quantum 
protocol serves as a superluminal communicator for transmitting any desired known classical information to any distant 
locations in the universe and that too instantaneously and simultaneously. The outcome of the protocol developed in section XI is not only far more 
important and useful but also very strange and shocking. This is so because in this protocol what information is to be sent to various participants,  
situated faraway from Alice and also from one another, is not known a priory. Here Alice can send any information generated 
long after Alice and all other participants have got separated from one another and also Alice can send any new information 
that will come into existence in future to all these distant participants, spacelike separated  from Alice 
at data transmission center. The other important application of this new quantum protocol for superluminally communicating any 
desired classical information in terms of classical bits, ordered in any desired sequence, is in achieving superluminal 
teleportation of quantum state from Alice to Bob as is seen above in section XII.

 \end{document}